\documentclass[reqno,12pt]{amsart}
\usepackage{amsmath, amssymb, amsthm, amsfonts} 
\usepackage{graphicx}
\usepackage[ruled,vlined]{algorithm2e}
\usepackage[a4paper,bindingoffset=0.5cm,left=2cm,right=2cm,top=2.5cm,bottom=2cm,footskip=.8cm]{geometry}
\usepackage{mathrsfs}

\newcommand{\bs}{\boldsymbol}

\newcommand{\vb}{\vspace{3.2mm}}
\renewcommand{\hat}{\widehat}

\newcommand{\vertiii}[1]{{\left\vert\kern-0.25ex\left\vert\kern-0.25ex\left\vert #1 \right\vert\kern-0.25ex\right\vert\kern-0.25ex\right\vert}}
\allowdisplaybreaks
\usepackage{hyperref}

\allowdisplaybreaks

\newtheorem{remark}{Remark}

\newtheorem{proposition}{Proposition}

\usepackage[utf8]{inputenc}
\usepackage{type1cm}         
\usepackage{multicol}        
\usepackage[bottom]{footmisc}

\usepackage{newtxtext}       %
\usepackage{newtxmath}       
\usepackage{multirow}

\usepackage{pgfplots}
\pgfplotsset{width=10cm,compat=1.9}
\usepgfplotslibrary{external}
\usepackage{subcaption}

\begin{document}

\title[Inference for dynamic Erd\H{o}s-R\'enyi random graphs under regime switching]{Inference for dynamic Erd\H{o}s-R\'enyi random graphs\\ under regime switching}

\author{Michel Mandjes
{\tiny and} Jiesen Wang}
	
\begin{abstract}
This paper examines a model involving two dynamic Erdős–Rényi random graphs that evolve in parallel, with edges in each graph alternating between being present and absent according to specified on- and off-time distributions.
A key feature of our setup is regime switching: the graph that is observed at any given moment depends on the state of an underlying background process, which is modeled as an alternating renewal process.
This modeling framework captures a common situation in various real-world applications, where the observed network is influenced by a (typically unobservable) background process. Such scenarios arise, for example, in economics, communication networks, and biological systems.

\noindent In our setup we only have access to aggregate quantities such as the number of active edges or the counts of specific subgraphs (such as stars or {complete graphs}) in the observed graph; importantly, we do not observe the mode. 
The objective is to estimate the on- and off-time distributions of the edges in each of the two dynamic Erdős-Rényi random graphs, as well as the distribution of time spent in each of the two modes. By employing parametric models for the on- and off-times and the background process, we develop a method of moments approach to estimate the relevant parameters. Experimental evaluations are conducted to demonstrate the effectiveness of the proposed method in recovering these parameters.

\vb

\noindent
{\sc Keywords.} Regime switching $\circ$  dynamic Erd\H{o}s-R\'enyi random graph $\circ$ method of moments $\circ$ subgraph counts

\vb

\noindent
{\sc Affiliations.} MM is with the Mathematical Institute, Leiden University, P.O. Box 9512,
2300 RA Leiden,
The Netherlands. He is also affiliated with Korteweg-de Vries Institute for Mathematics, University of Amsterdam, Amsterdam, The Netherlands; E{\sc urandom}, Eindhoven University of Technology, Eindhoven, The Netherlands; Amsterdam Business School, Faculty of Economics and Business, University of Amsterdam, Amsterdam, The Netherlands.


\noindent
JW is with the Korteweg-de Vries Institute for Mathematics, University of Amsterdam, Science Park 904, 1098 XH Amsterdam, The Netherlands.

\vb

\noindent
{\sc Acknowledgments.} 
This research was supported by the European Union’s Horizon 2020 research and innovation programme under the Marie Sklodowska-Curie grant agreement no.\ 945045, and by the NWO Gravitation project NETWORKS under grant agreement no.\ 024.002.003. \includegraphics[height=1em]{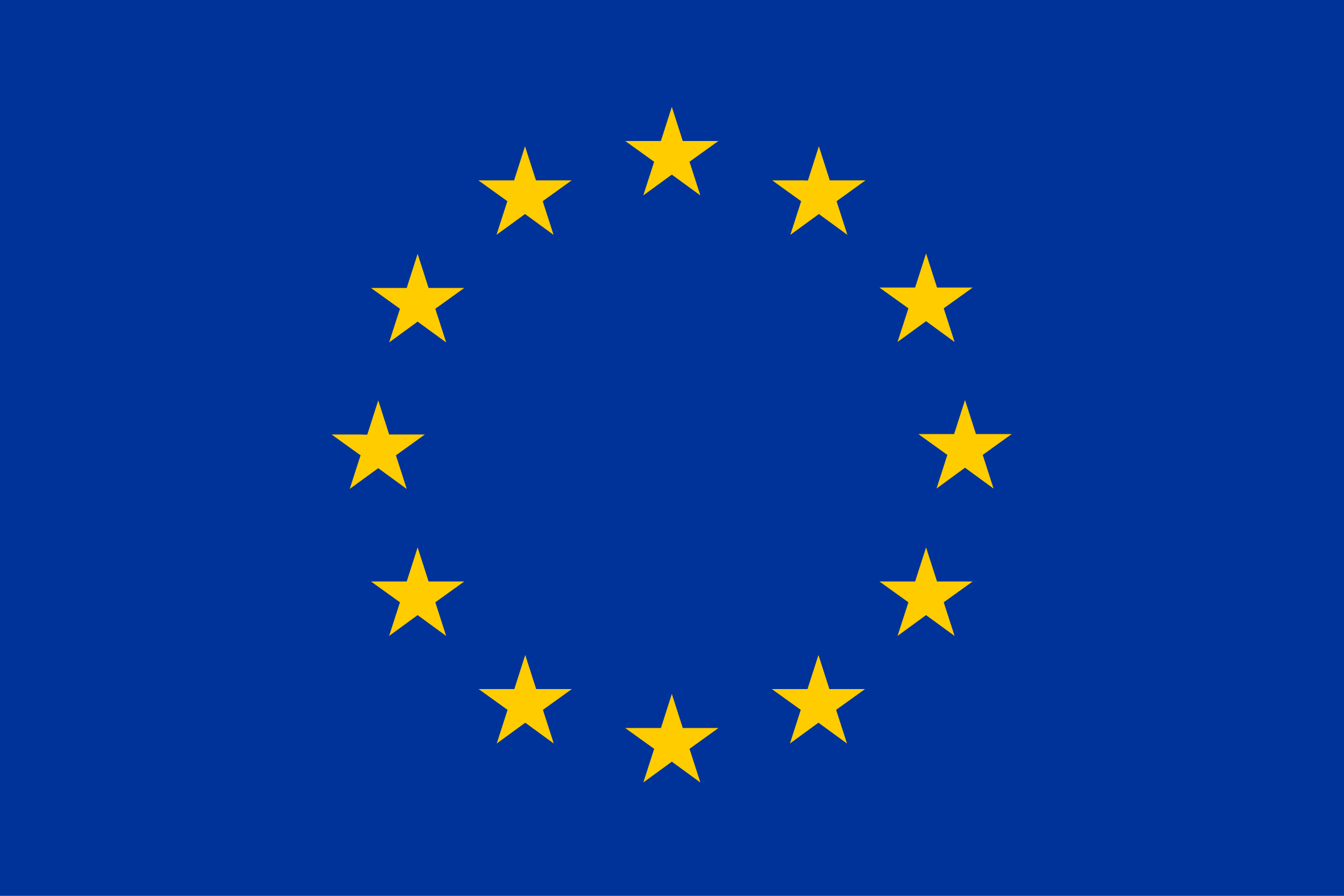} 
Date: {\it \today}.

\vb

\noindent
{\sc Email.} \url{m.r.h.mandjes@math.leidenuniv.nl}
and \url{j.wang2@uva.nl}.

\end{abstract}

\maketitle

\newpage 

\section{Introduction}

In the modeling of complex networks, various types of underlying stochastic mechanisms have been proposed. While traditional random graph models typically involve a single, {\it static} network realization, many real-world networks are of an intrinsically {\it dynamic} nature. This explains why in recent studies various dynamic versions of random graph models have been introduced, such as the dynamic counterpart of the conventional, static Erdős–Rényi random graph \cite{M2018}. 

In this paper we consider a model involving two dynamic Erdős-Rényi random graphs that run in parallel. The edges in each of these graphs alternate between being present and absent according to specific on- and off-time distributions. 
A distinguishing element of our setup is that there is {\it regime switching}, in that which of the two graphs is observed depends on the mode of a {\it background process} at that point in time. Throughout this paper, this background is modeled as an alternating renewal process. 
Assuming that all random variables involved belong to given parametric classes, the primary objective of this work is to develop an estimation procedure based on partial information derived from the observed process.



\vb

{\em Model description, observation scheme,  and objective}. We study a discrete-time model involving two independently evolving dynamic Erd\H{o}s-R\'enyi graphs, each having $N \in{\mathbb N}$ vertices. The stochastic dynamics of these graphs are defined as follows: in each graph, the $n := {N\choose 2}$ vertex pairs independently alternate between being present and absent. For graph $i$, with $i=1,2$, the on-times and off-times of each vertex pair are independent and identically distributed (i.i.d.), following the distributions of the generic random variables $X_i \in \mathbb{N}$ and $Y_i \in \mathbb{N}$. Additionally, the four sequences (the on-times and off-times for both graphs) are assumed to be independent.


At each discrete time point, which of the two dynamic Erdős–Rényi graph we observe is determined by the {\it mode} of a {\it background process}. This background process is of the regime switching type, in that the mode jumps between states 1 and 2 as an alternating renewal process, where the times spent in state $i$, for $i=1,2$, are distributed as the generic random variable $Z_i\in{\mathbb N}$. The resulting two sequences consist of i.i.d.\ random variables, and are assumed to be independent of each other as well as independent of the four sequences that we defined above. The mechanism underlying the background process is in the literature often referred to as `modulation'.


Throughout this study, we assume that the background process is {\it not observed}. Instead, we are just given, for a given set of subgraphs, the subgraph counts in the observed graph at each point in time (where we do not know which of the two graphs we are observing). For instance, it could be that we are given the time series pertaining to the number of edges and the number of triangles in the observed graph.  

In the setup that we consider, we assume that our system is {\it in stationarity}. This concretely entails that each of the edges in graph $i$ probabilistically behaves as a (static) Erd\H{o}s-R\'enyi random graph with $N$ vertices and ‘on-probability'
\[
\varrho_i \equiv \varrho_{X_i,Y_i}:= \frac{\mathbb{E} \, X_i}{\mathbb{E} \, X_i+\mathbb{E} \, Y_i},
\]
where we assume that the means $\mathbb{E}\,X_i, \, \mathbb{E}\,Y_i$ are finite, with $i=1,2$. Also, the regime process is in mode 1 with probability 
\[
\pi_1 \equiv \pi_{Z_1, Z_2} := \frac{\mathbb{E} \, Z_1}{\mathbb{E} \, Z_1+ \mathbb{E} \, Z_2} \,,
\]
and the system is in mode 2 with probability $\pi_2 := 1 - \pi_1$, where the means $\mathbb{E}\,Z_1, \, \mathbb{E}\,Z_2$ are assumed to be finite. 

Our objective is to {\it parametrically} estimate the on and off-time distributions in each of the two dynamic random graphs, as well as the distributions of $Z_1$ and $Z_2$, by observing the sequence of a specific subgraph counts over time. 
This means that we throughout assume that the six random variables that define our model belong to known, parametric families. Subgraph counts that are frequently used in the literature are the number of edges, wedges, or triangles.

A concrete example of our estimation problem is: in a setting in which we know that $X_i$ has a geometric distribution with parameter $p_i$, $Y_i$ a geometric distribution with parameter $q_i$ {and $Z_1,Z_2$ geometric distributions with parameters $p_0$ and $q_0$, respectively.}
{We} wish to estimate these six parameters based on the evolution of the number of triangles in the observed graph. 

\vb

{\em Subgraph counts}. As mentioned above, the input of our estimation procedure corresponds to a time series that records the evolution of specific subgraph counts in the observed graph. 
In the sequel we restrict ourselves to the following types of subgraphs:
\begin{itemize}
    \item[$\circ$] A {\it complete graph} $K_\ell$ on $\ell\in\{2,3,\ldots\}$ vertices is a simple undirected graph in which every vertex-pair is connected by a single edge. Thus $K_2$ is an edge and $K_3$ is a triangle. 
    \item[$\circ$] A {\it star} $S_\ell$ with $\ell\in\{2,3,\ldots\}$ is a tree with one `internal vertex' and $\ell-1$ `leaves'. Thus $S_2$ is an edge $S_2$ is an edge,  and $S_3$ is a path of length two, also known as a wedge in an undirected graph.
\end{itemize}

We proceed by introducing some convenient notation. In the first place, in self-evident notation, the following indicator function plays a crucial role throughout this paper:
\[\mathbf{1}^{(i)}_{j_1,j_2}(k)={\bs 1}\{\mbox{edge between vertices $j_1,j_2 $ in graph $i$ at time $k$ exists}\},\]
for vertices $j_1,j_2 \in \{1,\ldots, N\}$ with $j_1\not =j_2$, mode $i \in \{1,2\}$, and time $k \in \mathbb{N}$.
Denote the system mode at time $k\in{\mathbb N}$ by $M(k)$. The numbers of complete graphs or stars can be expressed in terms of the objects that we just introduced, as follows. Although the edge-count is a special case of both $K_\ell$ and $S_\ell$, we mention it separately, as it features prominently throughout this paper. 
\begin{itemize}
    \item[$\circ$] Aggregate number of edges at time $k$:
    \[\displaystyle A_N(k):= \sum_{i = 1}^2 \,\mathbf{1}_{\{M(k) = i\}} \,\sum_{j_1 = 1}^{{N\choose 2}-1} \sum_{j_2 = j_1+1}^{{N\choose 2}}\mathbf{1}^{(i)}_{j_1,j_2}(k) \,,\]
    \item[$\circ$] Aggregate number of $\ell$-complete graphs: \[\displaystyle K_{N,\,\ell}(k):= \sum_{i=1}^2
    \,\mathbf{1}_{\{M(k) = i\}} \,
    \sum_{\{[S]^\ell\}} \prod_{j_1, j_2 \in [S]^\ell} \mathbf{1}_{j_1,j_2}^{(i)}(k) \,,\]
    where $[S]^\ell$ is a subset of $\{1,\ldots,N\}$ with $\ell$ elements. Recalling that an edge is a $2$-complete graph, we have $A_N(k) = K_{N,\,2}(k)$.  
    \item[$\circ$] Aggregate number of $\ell$-stars:
    \[S_{N,\,\ell}(k) := \displaystyle \sum_{i=1}^2 \,\mathbf{1}_{\{M(k) = i\}} \,\sum_{\{[S]^\ell\}} \sum_{j_1 \in [S]^\ell} \prod_{j_2 \in [S]^\ell \backslash j_1}\mathbf{1}^{(i)}_{j_1,j_2}(k) \,.\]
    Recalling that an edge is also a $2$-star, we have $A_N(k) = S_{N,\,2}(k)$.
\end{itemize}

We focus on three fundamental structural patterns since they offer meaningful insights across a wide range of application domains.
\begin{itemize}
    \item {\it Edge.} An edge represents a direct connection between two nodes and serves as the basic unit of interaction in a network. In social networks, for instance, edges correspond to friendships or connections, where a greater number of edges typically indicates denser interactions, which is characteristic of tightly connected communities.
    \item {\it Complete subgraph.} The presence of complete subgraphs typically signifies tightly knit groups, such as close friend circles or collaborative teams. In social and communication networks, complete subgraphs are often associated with trust and mutual engagement. Platforms like Facebook and LinkedIn leverage this to suggest groups. Conversely, fraudulent or adversarial networks often lack such closed patterns, tending instead toward more random or sparse connectivity. In infrastructure networks, triangles indicate redundant paths, contributing to greater robustness and fault tolerance.
    \item {\it Star.} Stars are used to capture the potential for new connections and to identify central or influential nodes. In social networks, stars support recommendation systems such as Facebook’s `People You May Know' and LinkedIn’s connection suggestions. They also help uncover key individuals who serve as bridges between different communities. In financial transaction networks, a star-like structure (such as when User X transfers money to Y, and Y to Z, without a direct link between X and Z) may signal anomalous behavior and merit further investigation.
\end{itemize}

{\em Target objects}.
The following objects play a key role in our work: for $k \in \mathbb{N}, \ i \in \{1,2\}$, we define the six cumulative distribution functions by
\[
F_i(k):= \mathbb{P}(X_i \leqslant k), 
\qquad\qquad
G_i(k):= \mathbb{P}(Y_i \leqslant k),
\qquad\qquad
H_i(k):= \mathbb{P}(Z_i \leqslant k) \,.
\]
The corresponding probability mass functions, which can be mapped one-to-one to the cumulative distribution functions, are given by
\[
f_i(k):= \mathbb{P}(X_i = k), 
\qquad\qquad
g_i(k):= \mathbb{P}(Y_i = k), 
\qquad\qquad
h_i(k):= \mathbb{P}(Z_i = k) \,.
\]
It is throughout assumed that these six distributions belong to given parametric families, entailing that we seek to estimate a finite-dimensional parameter vector. 
As mentioned, we assume that the process is in stationarity, which in particular means that the system is in stationarity at time~1: at that epoch, each edge of dynamic random graph $i$ is on with probability $\varrho_i$ and off otherwise (for $i=1,2$), and the residual on-time $\bar X_i$, the residual off-time $\bar Y_i$, and the residual time $\bar Z_i$ in mode $i$ are characterized by the densities, for $k \in \mathbb{N}$ and  $i =1,2$,
\[
\bar f_i(k):= \frac{\sum_{\ell = k}^\infty f_i(\ell)}{\mathbb{E} \, X_i},
\qquad\qquad
\bar g_i(k):= \frac{\sum_{\ell = k}^\infty g_i(\ell)}{\mathbb{E} \, Y_i}
\qquad\qquad
\bar h_i(k):= \frac{\sum_{\ell = k}^\infty h_i(\ell)}{\mathbb{E} \, Z_i}
\,,
\]

respectively. 

We stress once more that we observe neither the evolution of individual edges nor the regime process. This concretely entails that we do not have direct access to realizations of the six distributions that we are interested in. Instead we have access to a time series corresponding to a vector of subgraph counts in the observed graph. It is noted that even in the specific (and probably most elementary) setting in which we observe the number of edges and all random variables are assumed to have geometric distributions, one does not have direct access to the underlying likelihood. Based on this observation, we deem that estimation by applying maximum likelihood is not a viable option.


\vb

{\em Test distributions}. Throughout this paper we work with various sorts of families of distributions, that we use to test our estimators. In the first place, we write $Z \sim \mathbb{G}(p)$ for the geometric random variable with ‘success probability’ $p$, for $p \in (0,1)$, when we mean
\[
\mathbb{P} (Z \geqslant k) = (1-p)^{k-1}, \qquad k = 1,2,\ldots.
\]
Second, we consider the more general class of Weibull distributions, covering tails that are lighter
and heavier than geometric. We write $Z \sim \mathbb{W}(\lambda, \alpha)$, with $\alpha, \lambda >0$, to denote 
\[
\mathbb{P} (Z \geqslant k) = e^{-\lambda(k-1)^\alpha}, \qquad k = 1,2,\ldots.
\]

\vb

{\em Contributions.} We study an inverse problem for a network where a hidden underlying process controls the observations. We assume the problem is parametric and adopt the method of moments for the estimation. Our contributions are summarized below.
\begin{itemize}
    \item[$\circ$] As mentioned, the observations correspond to subgraph counts. These could concern in principle any subgraph, but we explicitly work out the cases of the number of edges, complete graphs, and stars (or any combination of these). Setting up the moment equations requires the evaluation of a specific moment.  We derive closed-form expressions for specific types of (cross) moments, which form a crucial element in our approach. 
    \item[$\circ$] The performance of our estimation procedure is assessed through a series of experiments. Notably, our approach can deal with relatively high numbers of parameters; in one of the examples, we succeed in estimating as many as seven unknown parameters. Our experiments also shed light on the effect of the choice of the type of subgraph on the standard deviation of the estimator. 

\end{itemize}

{\it Organization.} {This paper is organized as follows. In Section~\ref{sec:literature} we review the literature on processes modulated by hidden underlying processes, focusing on applications and methodologies.
Then Section~\ref{sec:mom} introduces our approach based on the method of moments, where the estimation procedure relies on edge counts, complete subgraph counts, and star counts. Numerical results illustrating the performance of our estimator are presented Section~\ref{sec:NumericalExamples}. Section~\ref{sec:conclusion} concludes the paper with an account of future research directions.}
\section{Literature review}
\label{sec:literature}

In our model, observations are controlled by the underlying regime process. We are interested in estimating the parameters governing mode transitions and the graphs' dynamics. If the mode is observable, then the problem can be dealt with in the way discussed in \cite{MW23}. The crucial complication that we are facing now is that the mode is unobservable, which means that we do not know which of the two random graph processes each observation corresponds to.  

\medskip

{{\it Application areas.} Across a broad range of real-world applications we are facing the situation of observing a system governed by a non-observable modulating process, for instance in economics, communication networks, and biology. 
We proceed by providing a non-exhaustive account of papers in this domain. 
In \cite{Ou13} the focus is on fMRI images that reveal resting-state functional connectivity and its transition patterns in brain networks, providing insights into neural dynamics. Importantly, while the functional state of brain networks is unobservable, fMRI can measure blood-oxygen-level-dependent (BOLD) signals as an indirect indicator. In financial markets, stock returns are frequently modeled as being affected by the unobservable market state, which alternates between the so-called {\it bull} and {\it bear} markets. For additional applications in finance, see e.g.\ \cite{M14}.
In the context of communication networks, one often models  \cite{HL86, SV01} packet arrival processes as being controlled by a background process that switches between distinct states, each associated with a specific packet arrival rate (e.g.\  high, low, or no traffic). In various settings, one cannot observe the driving background process, but one {\it can} access performance measurements \cite{MvdM09}. In wireless networks, the unobservable background process could describe the atmospheric conditions. In biology, mRNA copy numbers are modulated by the active or inactive states of a gene; see e.g.\ the discussion in \cite{SRB12}. In mechanical systems, vibration signals collected from bearings reflect their health states, offering an indirect observation of the underlying condition \cite{OL05}.}

\medskip

{\it Methodologies.} In various studies methods have been devised for estimating the unobservable modulating process from observations of the process that it influences; see e.g.\ the survey \cite{AT} in the context of queues. An example is \cite{MM00}, in which the bull and bear market states, along with transitions between them, are modeled using a Markov switching framework, with parameters that are estimated by likelihood maximization.
Another approach is followed in \cite{HL86}, considering a setup in which packet arrival processes are modeled as Markov-modulated processes, with the parameters being inferred by matching characteristics of the observed process, such as the mean arrival rate, the variance-to-mean ratio, and the third moment of arrival counts. Similarly, \cite{GKMS19} examines a population process with independent sojourn times driven by a Markov-modulated Poisson process arrival mechanism. In queueing theory, this process corresponds to the M/M/$\infty$ queue in a random environment, as discussed in \cite{OP86}. These models are often interpreted as birth-death processes under modulation, with applications across various domains \cite{SRB12, AK15, OP86}. For such models, an Expectation-Maximization algorithm is frequently employed to obtain maximum likelihood estimates.
In communication networks, \cite{SV01} develops a methodology based on a hidden Markov model (HMM) so as to infer network states (e.g., congested or non-congested) by observing packet losses and delays in communication channels. Parameter estimation for the HMM is performed using the Expectation-Maximization algorithm. Similarly, \cite{OL05} utilizes an HMM trained in a supervised learning framework to analyze vibration signals from bearings, enabling health monitoring and fault prediction.
In the context of gene regulation, \cite{SRB12} assumes Gamma-distributed on and off times for genes, estimating parameters by matching the first three moments with their empirical counterparts.

\medskip

None of the frameworks described above covers our setting of dynamic random graphs under regime switching. Our problem is distinct due to the following features:
\begin{itemize}
\item We do not observe individual edge dynamics; instead, we only observe aggregate information, i.e. the total number of active edges at each time point. This aggregation renders many standard methods inapplicable.
\item The on and off-times of each edge can follow arbitrary distributions, not restricted to any particular parametric family.
\end{itemize}
To the best of our knowledge, no existing method can handle this combination of features. Our approach is based on the method of moments, incorporating moments that capture transitions between subsequent snapshots. Importantly, we estimate parameters for both the hidden underlying process and the observed process. 

\medskip

\section{Method of moments}
\label{sec:mom}

In this section, we propose an estimator for the unknown parameters based on the method of moments. The main idea of the method of moments is that we equate (i)~theoretical moment expressions (i.e., expressions for certain expectations, in terms of the parameters to be estimated), and (ii)~their empirical counterparts (i.e., estimates of these expectations, as obtained from observations). This leads to a number of equations, from which estimates of the unknown parameters can be derived.
Evidently, to use this methodology, it is essential that the number of moment equations is sufficiently large so as to be able to estimate the unknown parameters.

\medskip

Aiming at generating moment equations, 
the first natural candidates are the expected values pertaining to `single snapshots'. Indeed, for the mean number of edges, complete graphs and stars it is a straightforward exercise to verify that we have
\begin{align} \label{eq:mom1}
    & \mathbb{E} \, A_N(k) = {N \choose 2}\left(\pi_1 \,\varrho_1 + (1-\pi_1) \, \varrho_2 \right) \notag\\
    & \mathbb{E} \, K_{N,\ell}(k) = {N \choose \ell} \, \left(\pi_1 \, \varrho_1^{\ell \choose 2} + (1-\pi_1) \, \varrho_2^{\ell \choose 2} \right) \\
    &\mathbb{E} \, S_{N,\ell}(k) = \ell \, {N \choose \ell} \, \left(\pi_1 \, \varrho_1^{\ell-1} + (1-\pi_1) \, \varrho_2^{\ell-1} \right) \,, \notag
\end{align}
recalling that the process starts in stationarity. 
It is noted, however, that these expressions only contain the `equilibrium quantities' $\pi_1, \, \varrho_1,$ and $\varrho_2$, rendering our parameters unidentifiable. This issue we remedy by working with expectations involving {\it multiple} time epochs. 

\medskip

Throughout this paper, we work with `cross moments', i.e., the moments of quantities of the type  $A_N(k)A_N(k+d)$, $K_{N,\ell}(k)K_{N,\ell}(k+d)$, and $S_{N,\ell}(k)S_{N,\ell}(k+d)$, for $d \in {\mathbb N}$, on top of the `single snapshot moments' that we discussed above.  As will become clear below, a key challenge in calculating their expectations amounts to the counting of various patterns, to soundly reflect all relevant occurrences. In the remainder of this section, we primarily focus on the {\it lag-one} cross moments, i.e., the case that $d=1$. 

A first observation is that the probability that, in graph $i$, a specific edge is present at both time instant $k$ and time instant $k+1$ is
\[
\varrho_i\, (1-\bar f_i(1)),
\]
for $i=1,2.$
Here $\varrho_i$ is to be interpreted as the probability that in stationarity the edge in graph $i$ is present, and the factor $1-\bar f_i(1)$ represents the probability that it remains present. The probability that, in graph $i$, a specific edge is present at time $k$ and {\it another} edge at time $k+1$, is given by $\varrho_i^2$, with $i=1,2$, by virtue of each graph's edges evolving independently. Since these probabilities differ, in our computations of the cross moments, it is necessary to distinguish these cases.

\medskip

It is directly seen that the probability of observing different graphs at times $k$ and $k+1$ is given by \[r_{\not=}:=\pi_1 \bar h_1(1) + (1-\pi_1) \bar h_2(1).\] 
Similarly, the probability of observing graph $i$ at both time $k$ and time $k+1$ is $r_{i}:=\pi_i (1-\bar h_i(1))$, for $i=1,2$. Below we calculate the number of occurrences by distinguishing with respect to the number of {\it shared vertices}. In the sequel we denote by $a_m^{(\ell)}$ the number of two sets $\{i_1,\ldots, i_\ell\}$ and $\{j_1,\ldots, j_\ell\}$ that have $m$ elements in common. At any time $k$, there are ${N \choose \ell}$ ways to choose $\ell$ vertices, resulting in \[a_m^{(\ell)} = {N \choose \ell}{N-\ell \choose \ell-m}{\ell \choose m}\] for $m\in\{0,\ldots,\ell\}$. Also, it is readily verified that \[\sum_{m = 0}^\ell a_m^{(\ell)} = {N \choose \ell}^2 \,,\]
as desired.

\subsection{Edges}  Considering the case of $A_N(k)$, it is noted that, evidently, an edge is formed by selecting two vertices.
If two edges share no vertices or only one vertex, they are considered distinct. Conversely, if two edges share both vertices, they can be treated as being identical. Consequently, the cross moment $\mathbb{E} \, A_N(k)A_N(k+1)$ can be decomposed into three terms:
\begin{align} \label{eq:MoMA1}
    \mathbb{E} \, A_N(k)A_N(k+1)  = r_{\not=} \,{N \choose 2}^2 \, \varrho_1 \varrho_2 &+\, r_1 \, \left({a_2^{(2)}} \, \varrho_1 (1-\bar f_1(1)) + \left({a_0^{(2)}}+{a_1^{(2)}} \right)\varrho_1^2 \right)  \\
    &+\, r_2  \, \left({a_2^{(2)}} \, \varrho_2 (1-\bar f_2(1)) + \left({a_0^{(2)}+a_1^{(2)}} \right)\varrho_2^2 \right) \,. \notag
\end{align}
This decomposition can be understood as follows. The first term in the right-hand side of Equation~\eqref{eq:MoMA1} corresponds to the scenario where at time $k$ and time $k+1$ two different graphs are observed. Since each graph contains ${N \choose 2}$ edges, this can occur in ${N \choose 2}\times{N \choose 2} $ possible ways.  The probability of observing a specific edge at both times $k$ as well as $k+1$ is given by $\varrho_1 \, \varrho_2$, as a consequence of the assumed independence and stationarity. The second term in the right-hand side of Equation~\eqref{eq:MoMA1}  corresponds to the scenario where graph 1 is observed at time $k$ as well as at time $k+1$. The probability of a specific edge being present at both time instances is $\varrho_1(1-\bar f_1(1))$ if the edge is the same (hence the weight $a_2^{(2)}$), and otherwise this probability is $\varrho_1^2$ (with the weight $a_1^{(2)}$). This reasoning also applies to the third term in~\eqref{eq:MoMA1}, covering the scenario where it is graph 2 that is observed at the time instances  $k$ and $k+1$. 

\subsection{Complete graphs} In the case of the complete graph $K_{\ell}$, forming it requires  $\ell$ vertices to be connected.
Given $m$ common vertices,  the common edges are the edges among these $m$ vertices. With $b_{\ell}:= {N \choose \ell}$ denoting the total possible edges for $\ell$ vertices, the number of common edges is $b_m$, leading to
\begin{align} \label{eq:MoMK1}
    \mathbb{E} \, K_{N, \ell}(k)K_{N, \ell}(k+1) 
    = r_{\not=} \,{N \choose \ell}^2 \varrho_1^{b_{\ell}} \varrho_2^{b_{\ell}}  
    &+ \, r_1\, \sum_{m=0}^{\ell} \left({a_m^{(\ell)}} \, \varrho_1^{b_{\ell}} \,(1-\bar{f}_1(1))^{b_m} \,  \varrho_1^{b_{\ell}-b_m} \right) \\ &+ \, r_2 \, \sum_{m=0}^{\ell} \left( {a_m^{(\ell)}} \, \varrho_2^{b_{\ell}} \, (1-\bar{f}_2(1))^{b_m} \, \varrho_2^{b_{\ell}-b_m} \right) \notag \,. \notag
\end{align}   
The interpretation of this expression is similar to the one underlying ~\eqref{eq:MoMA1}. More specifically, the first term in the right-hand side of Equation~\eqref{eq:MoMK1} corresponds to the scenario that at times $k$ and $k+1$ two different graphs are observed. Now there are ${N \choose \ell}\times {N \choose \ell}$ possible configurations, and the probability of observing a graph of the type $K_{\ell}$ at both times is given by $\varrho_1^{b_{\ell}} \, \varrho_2^{b_{\ell}}$.  The second term in the right-hand side of Equation \eqref{eq:MoMK1} corresponds to the scenario that graph~1 is observed at time $k$ as well as at time $k+1$. When two complete subgraphs have $m$ common vertices, then they share $b_m$ common edges. A complete subgraph is observed at time $k$ with probability $\varrho_1^{b_{\ell}}$. Given it is observed, another complete subgraph that shares $m$ vertices with it is observed with probability $(1- \bar f_1(1))^{b_m}\varrho_1^{b_{\ell}-b_m}$. This reasoning also applies to the third term in~\eqref{eq:MoMA1}, which corresponds to the scenario where graph~2 is observed at both instances $k$ and $k+1$.

\subsection{Stars} In the case of the star graph $S_{\ell}$, forming it concerns $\ell$ vertices.
Given $m$ common vertices, the number of common edges is more involved than in the case of the complete graph, due to the star's intrinsic asymmetry (i.e., one has to distinguish between leaf vertices and the center vertex). In Proposition \ref{prop:1} we provide, for two stars $S_{\ell}$ and $\Tilde{S}_{\ell}$, the number of common edges when they share $m$ vertices.
Denote by ${\mathscr C}_1$ the scenario that at most one center vertex is a common vertex, by ${\mathscr C}_2$ the scenario that two center vertices are common vertices but that these vertices are different, and by ${\mathscr C}_3$ the scenario that the two center vertices are the same.

\begin{proposition} \label{prop:1}
Let the two stars $S_{\ell}$ and $\Tilde{S}_{\ell}$ share $m$ vertices. Then Table \ref{tab:1} presents, for the scenarios ${\mathscr C}_1$, ${\mathscr C}_2$, and ${\mathscr C}_3$, the number of common edges and the number of occurrences.
\begin{table}[ht]
\caption{Number of common edges in different cases.}
\centering
  \begin{tabular}{||c c c||} 
    \hline
      & \mbox{\rm number of common edges} & \mbox{\rm number of cases}    \\ [0.5ex] 
     \hline\hline
     \mbox{\rm ${\mathscr C}_1$} & $0$ & $(\ell+m)(\ell -m)$\\ 
     \hline
     \mbox{\rm ${\mathscr C}_2$} & $1$ & $m(m-1)$\\ 
     \hline
     \mbox{\rm ${\mathscr C}_3$} & $m-1$ & $m$   \\
     \hline
     \end{tabular}
  \label{tab:1}
\end{table}
\end{proposition}
\begin{proof}
    We analyze the number of common edges in three situations.
\begin{itemize}
\item[      $\circ$] If the centers of  $S_{\ell}$ and $\Tilde{S}_{\ell}$ are not among the common vertices, there are no common edges. In this case, each star has $(\ell-m)$ possible choices for its center, resulting in $(\ell-m)^2$ combinations. If the center of one star is a common vertex while the center of the other star is not, the first star has $m$ choices for its center, and the second star has $\ell - m$ choices. Since either star can take the role of the first, there are $2m(\ell-m)$ possibilities. Adding up the number of cases corresponding to these two scenarios, the total number of configurations corresponding to zero common edges is $(\ell+m)(\ell-m)$.

 \item[    $\circ$] If the centers of $S_{\ell}$ and $\Tilde{S}_{\ell}$ are different, but are both among the common vertices, then there is one sharing edge, that is the edge linking the two centers. In this situation, the two center vertices are from the $m$ vertices, so that there are $m(m-1)$ possibilities.

 \item[    $\circ$] If the center of  ${S}_{\ell}$ is also the center of $\Tilde{S}_{\ell}$, then the two subgraphs share $m-1$ common edges. As the common center can be any of the $m$ vertices, there are $m$ possibilities.
 \end{itemize}
 This concludes the proof. 
\end{proof}

\begin{figure}
    \centering    {\includegraphics[width=0.8\linewidth]{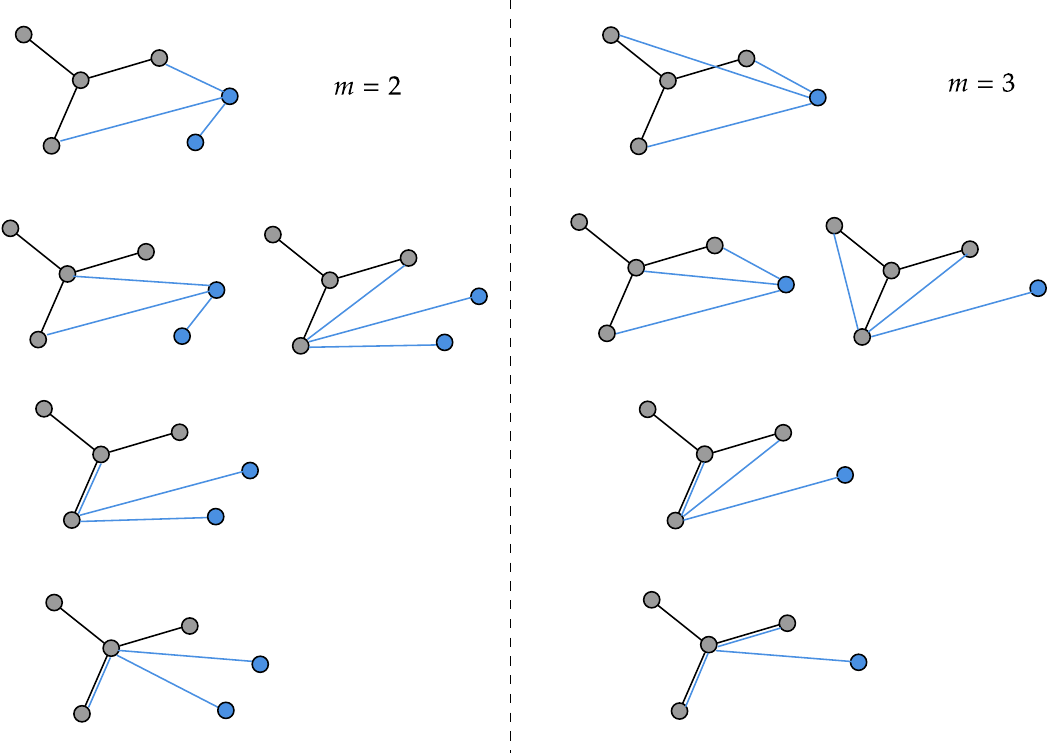}}
    \caption{Illustration of common edges for $S_{4}$ and $\Tilde{S}_{4}$ when $m = 2$ (left panels) and $m=3$ (right panels), where the grey edges and nodes represent ${S}_4$, and the blue edges and nodes represent $\Tilde{S}_4$.}
    \label{fig:StarCommon}
\end{figure}

We illustrate the result in Table \ref{tab:1} for the two star graphs $S_{4}$ and $\Tilde{S}_{4}$. In case $m = 0$ or $4$, it is directly seen that there are $0$ and $4$ common edges (as they are the same subgraph). So we can restrict ourselves, in Figure \ref{fig:StarCommon}, to $m = 2$ (left panels) and $3$ (right panels). In the first line, where the two center vertices are not among the $m$ common vertices, the two stars have no sharing edges. In the second line only either the center of ${S}_4$ (left) or the center of $\Tilde{S}_4$ (right) is a common vertex, and the other center is not. The two stars have no sharing edges in this case. In the third line the two center vertices are different common vertices, and there is one sharing edge which is the one linking the two center vertices. In the fourth line the center vertices of $S_4$ and $\Tilde{S}_4$ are the same, and there are $m-1$ sharing edges.

As mentioned before, a wedge can be regarded as the star $S_3$. According to Table \ref{tab:1}, we have
\begin{itemize}
    \item[    $\circ$] if  $m=0,1$, then there are 0 common edges.
    \item[    $\circ$] if $m=2$, then there are 5 cases with 0 common edges, and 4 cases with 1 common edge.
    \item[    $\circ$]  if $m = 3$, then there are 3 cases with 2 common edges, and 6 cases with 1 common edge.
\end{itemize}
These results are consistent with our findings in \cite{MW23}.

\medskip

It follows from Proposition \ref{prop:1} that the following decomposition applies:
\begin{align} \label{eq:MoMS1}
    &\mathbb{E} \, S_{N, \ell}(k)S_{N, \ell}(k+1)= r_{\not=} \left(\ell \,{N \choose \ell} \right)^2 \, \varrho_1^{\ell-1} \varrho_2^{\ell-1}  \\ &+r_1\, \sum_{m=0}^{\ell} {a_m^{(\ell)}} \left((\ell -m)(\ell +m) \varrho_1^{2(\ell-1)} + m(m-1) \varrho_1^{2\ell-3} (1-\bar f_1(1)) + m  \varrho_1^{2\ell -m-1}(1-\bar f_1)^{m-1} \right) \notag \\
    &+ r_2 \, \sum_{m=0}^{\ell} {a_m^{(\ell)}} \left((\ell -m)(\ell +m) \varrho_2^{2(\ell-1)} + m(m-1) \varrho_2^{2\ell-3} (1-\bar f_2(1)) + m  \varrho_2^{2\ell -m-1}(1-\bar f_2)^{m-1} \right) \,. \notag
\end{align}
Similarly to how we have interpreted Equations \eqref{eq:MoMA1} and \eqref{eq:MoMK1}, each of the individual terms of the right-hand side of Equation \eqref{eq:MoMS1} corresponds to a specific scenario. The first term reflects the scenario where two different graphs are observed at times $k$ and $k+1$. Forming a star  $S_\ell$ requires $\ell$ vertices, and any of these $\ell$ can be the central vertex, entailing that there are $\ell\, {N \choose \ell} $ ways to form a star in any of the two graphs. The star has $\ell-1$ edges between the central vertex and the leaf vertices, so the probability of observing the star is $\varrho_1^{\ell-1}$. The  two graphs evolve independently, so that we obtain the contribution \[\ell\,{N \choose \ell} \,\varrho_1^{\ell-1} \, \ell{N \choose
\ell} \, \varrho_2^{\ell-1}.\] The second term covers the scenario where graph 1 is observed at time $k$ and $k+1$. With the results given in Proposition \ref{prop:1}, we have \[a_m^{(\ell)} (\ell +m)(\ell -m)\] cases in which the two graphs at $k$ and $k+1$ sharing $m$ vertices have no edges in common. The probability of observing the two stars is $\varrho_1^{\ell-1} \, \varrho_1^{\ell-1}$. Again appealing to Proposition \ref{prop:1}, there are $a_m^{(\ell)} \, m(m-1)$ cases corresponding to the scenario where the two graphs at $k$ and $k+1$ sharing $m$ vertices have $1$ edge in common. The probability of observing two stars in this scenario is $\varrho_1^{\ell-1} \, \varrho_1^{\ell-2} \, (1-\bar f_1(1))$. There are $a_m^{(\ell)} m$ cases corresponding to the scenario where two stars sharing $m$ vertices have $m-1$ edges in common. The probability of observing these two stars is $\varrho_1^{\ell-1} \, \varrho_1^{\ell -m} \, (1-\bar f_1(1))^{m-1}$. Adding up these three parts, and multiplying by $r_1$, results in the second term in the right-hand side of \eqref{eq:MoMS1}. The third term in Equation \eqref{eq:MoMS1}, corresponding to the scenario where graph 2 is observed at $k$ and $k+1$, is derived along the same lines. 

\subsection{Estimator}\label{ssec:est}
Now that we have expressions for the single snapshot moments and the cross moments, we can propose our estimator. 
Suppose for the moment that there is a single unknown parameter for each of the random variables involved, so that there is a total of six unknown parameters that we would like to estimate. Theoretically, we require six independent moments involving the six unknown parameters, where in principle any value of $\ell$ could be used. However, smaller values of $\ell$ (a)~lead to simpler expressions for the moments, and (b)~autocorrelations that are more clearly separated from 0. Therefore, in our experimental evaluation, we focused on $\ell \in {2, 3, 4}$.

First observe that the equilibrium quantities $\pi_1, \, \varrho_1, \, \varrho_2$ are functions of the parameters pertaining to the random vectors $(X_1, Y_1)$, $(X_2, Y_2)$, and $(Z_1, Z_2)$, respectively. In our approach, we first solve for the three `equilibrium quantities' $\pi_1, \, \varrho_1, \, \varrho_2$ using observations of the time series $A_N(k)$, $K_{N, \ell}(k)$, and $S_{N, \ell}(k)$, based on the single snapshot moments. The second step is to estimate the six parameters by combining results from the first step with the cross-moments corresponding to $A_N(k)A_N(k+1)$, $K_{N, \ell}(k)K_{N, \ell}(k+1)$, and $S_{N, \ell}(k)S_{N, \ell}(k+1)$.

As mentioned, the moment equations are calculated in two steps. Let $T\in{\mathbb N}$ throughout denote the number of observations. 
\begin{itemize}
    \item[$\circ$] {\it Step 1}: First define
    \[
    \mathcal{A}_{T,0}:= \frac{1}{T} \sum_{k = 1}^T \, A_N(k), \quad \mathcal{K}_{T,\ell,0} := \frac{1}{T} \sum_{k = 1}^T \, K_{N, \ell}(k), \quad \mathcal{S}_{T,\ell,0} := \frac{1}{T} \sum_{k = 1}^T \, S_{N, \ell}(k) \,,
    \]
    which are evident estimators of the expectations $\mathbb{E} \, A_N(k), \,  \mathbb{E} \, K_{N, \ell}(k)$, and $\mathbb{E} \, S_{N, \ell}(k)$, respectively.
Then the estimators $\hat \pi_1$, $\hat \varrho_1$, and $\hat \varrho_2$ (of the probabilities $\pi_1$, $\varrho_1$, and $\varrho_2$, respectively) are obtained by solving the three equations
\begin{equation} \label{eq:MoMStep1}
    \mathbb{E} \, A_N(k) = \mathcal{A}_{T,0}, \quad  \mathbb{E} \, K_{N, \ell}(k) = \mathcal{K}_{T,\ell,0}, \quad \mathbb{E} \, S_{N, \ell}(k) = \mathcal{S}_{T,\ell,0} \,,
\end{equation}
where $\mathbb{E} \, A_N(k)$, $\mathbb{E} \, K_{N, \ell}(k)$, and $\mathbb{E} \, S_{N, \ell}(k)$ are as given in \eqref{eq:mom1}. The calculation in this step is the same for any distribution, i.e., it does not involve the specific parametric form of the six distributions. 

\item[$\circ$] {\it Step 2}: Then we define
\begin{align*}
    & \mathcal{A}_{T,1}:= \frac{1}{T-1} \sum_{k = 1}^{T-1} \, A_N(k)A_N(k+1)  ,\\
    & \mathcal{K}_{T,\ell,1} := \frac{1}{T-1} \sum_{k = 1}^{T-1} \, K_{N, \ell}(k)K_{N, \ell}(k+1) ,\\
    & \mathcal{S}_{T,\ell,1} := \frac{1}{T-1} \sum_{k = 1}^{T-1} \, S_{N, \ell}(k)S_{N, \ell}(k+1) \,,
\end{align*}
being estimators of $\mathbb{E} \, A_N(k)A_N(k+1)$, $\mathbb{E} \, K_{N, \ell}(k)K_{N, \ell}(k+1)$, and $\mathbb{E} \, S_{N, \ell}(k)S_{N, \ell}(k+1)$, respectively. 
It follows that the moment equations are 
\begin{align} \label{eq:MoM2}
    \mathcal{A}_{T,1} &= \mathbb{E} \, A_N(k)A_N(k+1) \notag , \\
    \mathcal{K}_{T,\ell,1} &= \mathbb{E} \, K_{N, \ell}(k)K_{N, \ell}(k+1) ,\\
    \mathcal{S}_{T,\ell,1} &= \mathbb{E} \, S_{N, \ell}(k)S_{N, \ell}(k+1) ,\notag 
\end{align}
where $\mathbb{E} \, A_N(k)A_N(k+1), \,  \mathbb{E} \, K_{N, \ell}(k)K_{N, \ell}(k+1)$, and $\mathbb{E} \, S_{N, \ell}(k)S_{N, \ell}(k+1)$ are as given in \eqref{eq:MoMA1}, \eqref{eq:MoMK1}, and \eqref{eq:MoMS1}.
It is readily verified from its definition that $\bar h_1(1)$ can be written as $1/\mathbb{E} Z_1$; similarly, $\bar h_2(1), \bar g_1(1),\bar g_2(1), \bar f_1(1), \bar f_2(1)$ can be written as the reciprocal of the corresponding expected values. Moreover, $\mathbb{E} Y_1, \, \mathbb{E} Y_2$ and $\mathbb{E} Z_2$ can be expressed by $\mathbb{E} X_1, \, \mathbb{E} X_2$ and $\mathbb{E} Z_1$ and the stationary distributions, i.e.,
\begin{equation} \label{eq:par1}
    \mathbb{E}\, Y_1 = \Theta(\varrho_1, \mathbb{E}\, X_1) \qquad \mathbb{E}\, Y_2 =  \Theta(\varrho_2, \mathbb{E} \, X_2) \qquad  \mathbb{E}\, Z_2 = \Theta(\pi_1, \mathbb{E} \,Z_1)
\end{equation}
where $\Theta(x,y) := (1-x)y/x$. 
In Step 1, we estimated the three stationary probabilities $\hat \pi_1$, $\hat \varrho_1$, and $\hat \varrho_2$ , so now we only need to determine $\mathbb{E} \, X_1, \, \mathbb{E} \,X_2$, and $\mathbb{E}\, Z_1$. 
\end{itemize}
Calculating the moment equations in two steps, rather than solving the six equations simultaneously, drastically reduces the computational complexity.  In summary, we first solve for $\hat\pi_1, \, \hat\varrho_1, \, \hat\varrho_2$ using \eqref{eq:MoMStep1}, then ${\mathbb{E} \, Z_1}, \, {\mathbb{E} \, X_1}, \, {\mathbb{E} \, X_2}$ via \eqref{eq:MoM2} and \eqref{eq:par1}. The values of ${\mathbb{E} \, Z_2}, \, {\mathbb{E} \, Y_1}, \, {\mathbb{E} \, Y_2}$ are subsequently inferred. Finally, the estimators of the six parameters are then derived from these estimates.

\begin{remark}\label{R1}{\em 
    Besides the moments mentioned above, we can also consider the observations at time instants $k$ and $k+2$. In this remark we consider, as an example, the case of the number of edges. Let
\[
\mathcal{A}_{T,2}:= \frac{1}{T-2} \sum_{k = 1}^{T-2} \, A_N(k)A_N(k+2) 
\]
be an evident estimator of the `lag-two cross moment'. To simplify the notation, define $P_{ij}(d):= \mathbb{P}(M(k+d) = j\,|\, M(k) = i)$. Then,
\begin{align*}
    & P_{11}(2) = (1-\bar h_1(1))^2 + \bar h_1(1) \bar h_2(1) \qquad P_{12}(2) = \bar h_1(1) (1-\bar h_2(1)) \\
    & P_{22}(2) = (1-\bar h_2(1))^2 + \bar h_2(1) \bar h_1(1)   \qquad P_{21}(2) =\bar h_2(1) (1-\bar h_1(1)) \,.
\end{align*}
The probability of observing different graphs at $k$ and $k+2$ is $\pi_1 P_{12}(2) + (1-\pi_1) P_{21}(2)$. If the same edge in graph 1 is observed at times $k$ and $k+2$,  
\begin{itemize}
    \item[$\circ$] with probability $1-\bar f_1(1) - \bar f_1(2)$, the edge remains on at times $k, k+1,$ and $k+2$;
    \item[$\circ$] with probability  $\bar f_1(1) g_1(1)$, the edge is on at time $k$, off at $k+1$, and on again at $k+2$.
\end{itemize}
Thus, the probability of observing an edge at time epochs $k$ and $k+2$ is \[\varrho_1 (1-\bar f_1(1) - \bar f_1(2) + \bar f_1(1)  g_1(1)).\] If two different edges are observed at $k$ and $k+2$, then the probability is $\varrho_1^2$. A similar reasoning applies to graph 2. Thus,
\begin{align} \label{eq:MoMA3}
    \mathbb{E} \, A_N(k)A_N(k+2) &= \left(\pi_1 P_{12}(2) + (1-\pi_1) P_{21}(2)\right) {N \choose 2}^2 \, \varrho_1 \varrho_2\ \\ \notag 
    &+ \pi_1 P_{11}(2) \, \left(n{a_2^{(\ell)}} \, \varrho_1 (1-\bar f_1(1) - \bar f_1(2) + \bar f_1(1)  g_1(1)) + \left({a_0^{(\ell)}}+{a_1^{(\ell)}} \right)\varrho_1^2 \right) \notag \\
    &+ (1-\pi_1) P_{22}(2) \, \left({a_2^{(\ell)}} \, \varrho_2 (1-\bar f_2(1) - \bar f_2(2) + \bar f_2(1)  g_2(1)) + \left({a_0^{(\ell)}}+{a_1^{(\ell)}} \right)\varrho_2^2 \right) \,. \notag
\end{align}

Theoretically, a prerequisite for the method of moments to work is that the number of moment equations equals the number of unknowns. However, with the time between the two observations increasing, the dependence between the observations becomes weaker. This makes, as $d$ grows, $\mathbb{E} \, A_N(k)A_N(k+d)$ `increasingly similar' to $\mathbb{E} \, A_N(k)\cdot \mathbb{E} \,A_N(k+d)$, so that the corresponding moment equation carries less and less information.  }
\end{remark}

\section{Numerical examples} 
\label{sec:NumericalExamples}
We illustrate the performance of our estimation approach through a series of 
experimental evaluations. Specifically, in these experiments we let $Z_1 \sim \mathbb{G}(p_0), \, Z_2 \sim \mathbb{G}(q_0)$ with $p_0 = 0.3, \, q_0 = 0.6$, and we consider the following three instances:
\begin{itemize}
       \item[$\circ$] Case I: $X_1 \sim \mathbb{W}(1.5, \alpha), \, Y_1 \sim \mathbb{G}(q_1), \, X_2 \sim \mathbb{W}(1.5, \beta),\, Y_2 \sim \mathbb{G}(q_2)$,  with $\alpha = 0.5, \, q_1 = 0.4, \, \beta = 0.3, \, q_2 = 0.8$ {to be estimated}. This makes six parameters, which we estimate using the moment equations based on $\mathcal{A}_{T,0}$, $\mathcal{K}_{T,3,0}$, $\mathcal{S}_{T,3,0}$ and $\mathcal{A}_{T,1}$, $\mathcal{K}_{T,3,1}$, $\mathcal{S}_{T,3,1}$.
    \item[$\circ$] Case II: $X_1 \sim \mathbb{W}(1.5, \alpha), \, Y_1 \sim \mathbb{G}(q_1), \, X_2 \sim \mathbb{W}(1.5, \beta),\, Y_2 \sim \mathbb{G}(q_2)$, with the same {unknown} parameters as in Case I. This makes again six parameters, but now we use the moment equations based on  $\mathcal{A}_{T,0}$, 
    $\mathcal{S}_{T,3,0}$, $\mathcal{S}_{T,4,0}$ and $\mathcal{A}_{T,1}$, $\mathcal{S}_{T,3,1}$, $\mathcal{S}_{T,4,1}$;
    \item[$\circ$] Case III: $X_1 \sim \mathbb{W}(\lambda, \alpha), \, Y_1 \sim \mathbb{G}(q_1), \, X_2 \sim \mathbb{W}(1.5, \beta),\, Y_2 \sim \mathbb{G}(q_2)$,  where, in addition to the six {unknown} parameters that we came across in Cases I and II, now we also wish to estimate
    $\lambda = 1.5$.  This means that there are seven parameters, which we estimate by using the moment equations based on $\mathcal{A}_{T,0}, \, \mathcal{K}_{T,3,0}, \, \mathcal{S}_{T,3,0}$ and $\mathcal{A}_{T,1}, \,\mathcal{A}_{T,2}, \, \mathcal{K}_{T,3,1}, \, \mathcal{S}_{T,3,1}$ (where for $\mathcal{A}_{T,2}$ we use Remark \ref{R1}).
\end{itemize}
For Cases I and II, in which there is one unknown parameter for each of the six random variables involved, we outlined in Section \ref{ssec:est} a procedure to estimate the six means ${\mathbb{E} \, Z_1}$, ${\mathbb{E} \, Z_2}$, ${\mathbb{E} \, X_1}$, ${\mathbb{E} \, X_2}$, ${\mathbb{E} \, Y_1}$, ${\mathbb{E} \, Y_2}$ from the six moment equations pertaining to single-snapshot moments and lag-one mixed moments.

It is an important observation that if we change the class of distributions of some of the random variables, the estimation procedure remains essentially the same. For instance compare the case that $X_1 \sim \mathbb{P}{\rm ar}(1, \alpha)$ with the case that $X_1 \sim \mathbb{G}(p_1)$. Indeed, in both cases we first estimate ${\mathbb{E} \,X_1}$ (besides the other five expectations), leading to an estimator that we call $\hat\mu.$ Then in case $X_1 \sim \mathbb{P}{\rm ar}(1, \alpha)$ we have $\hat\alpha = \zeta^{-1}(\hat\mu)$ with $\zeta(x) = \sum_{k = 1}^\infty k^{-x}$ denoting the Riemann zeta function, whereas in case $X_1 \sim \mathbb{G}(p_1)$ we have  $\hat p_1 = 1/{\hat\mu}$; here it has been used that for $X_1 \sim \mathbb{P}{\rm ar}(1, \alpha)$ we have \[{\mathbb{E} \,X_1}=
\sum_{k=1}^\infty {\mathbb P}(X_1\geqslant k) = \sum_{k = 1}^\infty \frac{1}{k^{\alpha}}=
\zeta(\alpha).\]

\medskip

We proceed by briefly discussing the motivation behind choosing these three instances. 
We chose in Cases I and II the same parameter setting, but the estimation is performed with different (lag-one) cross moments. The objective is to assess the effect of the subgraphs on the performance of the estimator. In Case III, we have seven (rather than six) unknown parameters; we include the lag-two cross moment ${\mathbb E}[A_N(k)A_N(k+2)]$, besides three single snapshot moments and three lag-one cross moments.

\begin{figure}
    \centering
    \subcaptionbox{$\hat{p}_0^{(L)} = 0.3000  \, (0.0020)$}
    {\includegraphics[width=0.49\linewidth]{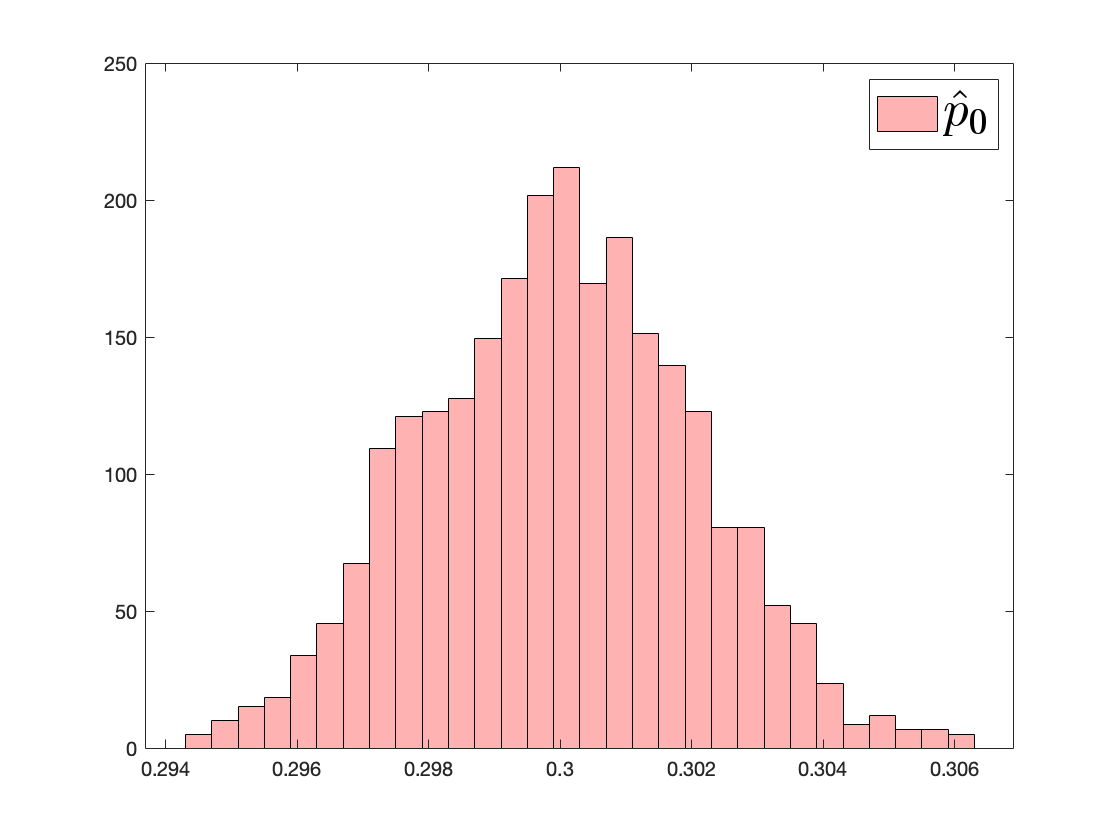}}
    \subcaptionbox{${\hat{q}_0^{(L)}} = 0.6002 \, (0.0030)$}
    {\includegraphics[width=0.49\linewidth]{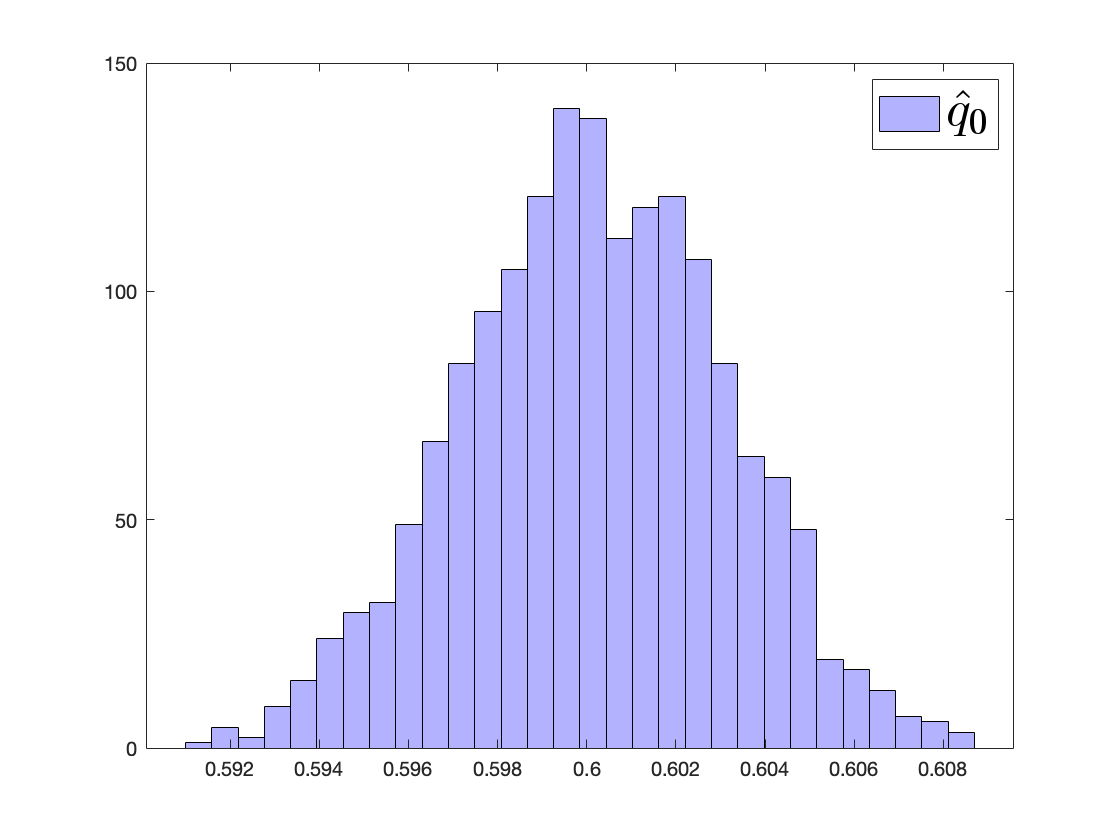}}\\
    \subcaptionbox{$\hat{\alpha}^{(L)} = 0.5003 \, (0.0127)$}
    {\includegraphics[width=0.49\linewidth]{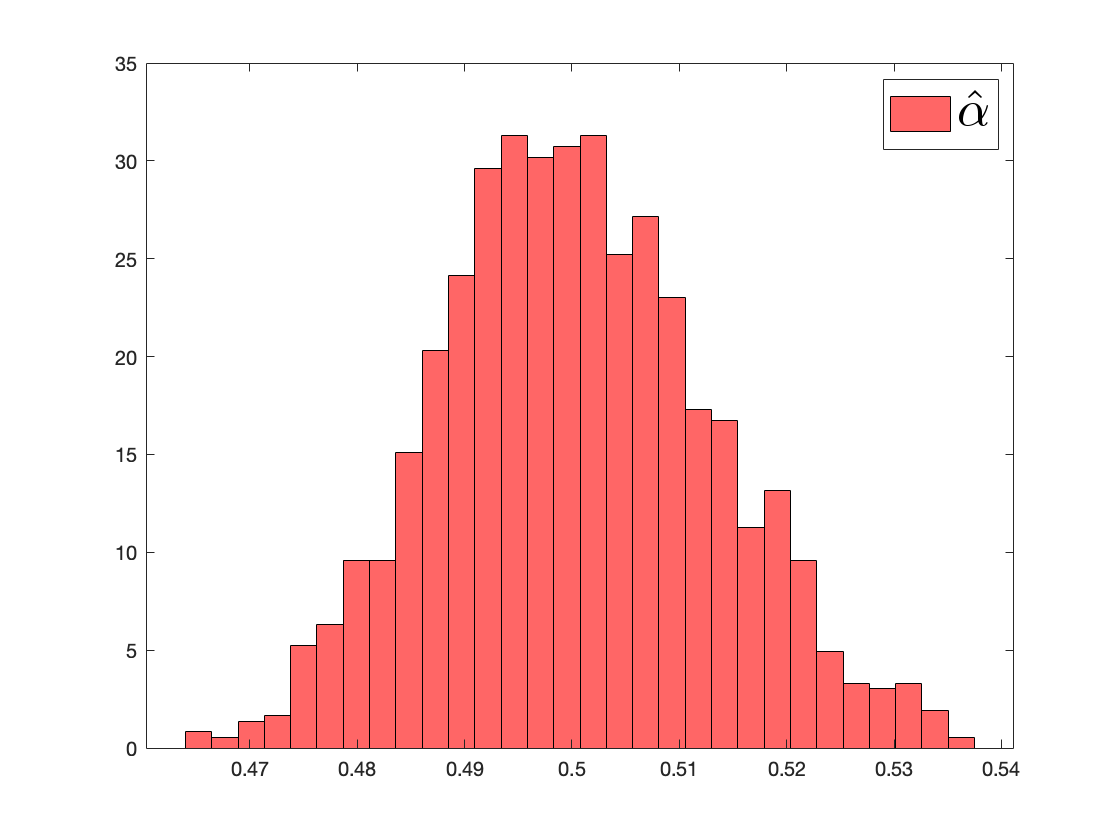}} 
    \hfill
    \subcaptionbox{$\hat{q}_1^{(L)} = 0.4000 \, (0.0059)$}
    {\includegraphics[width=0.49\linewidth]{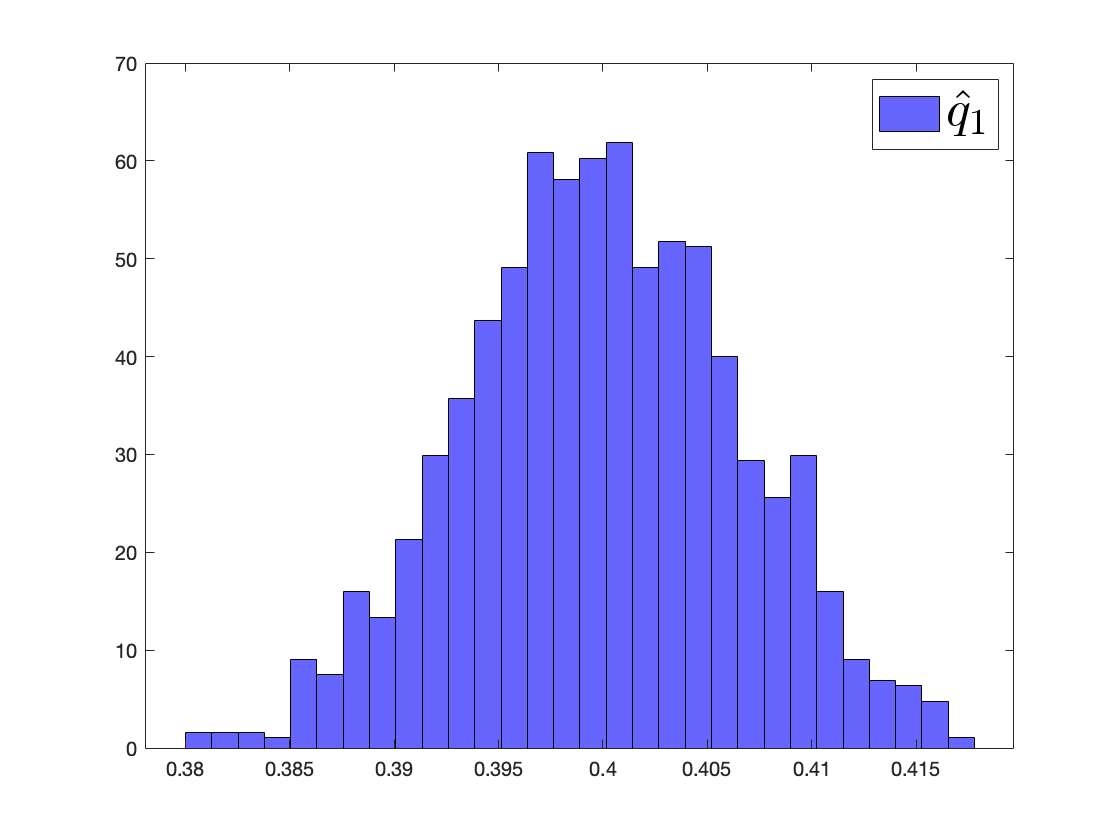}} \\
    \subcaptionbox{$\hat{\beta}^{(L)} = 0.3000 \, (0.0044)$}
    {\includegraphics[width=0.49\linewidth]{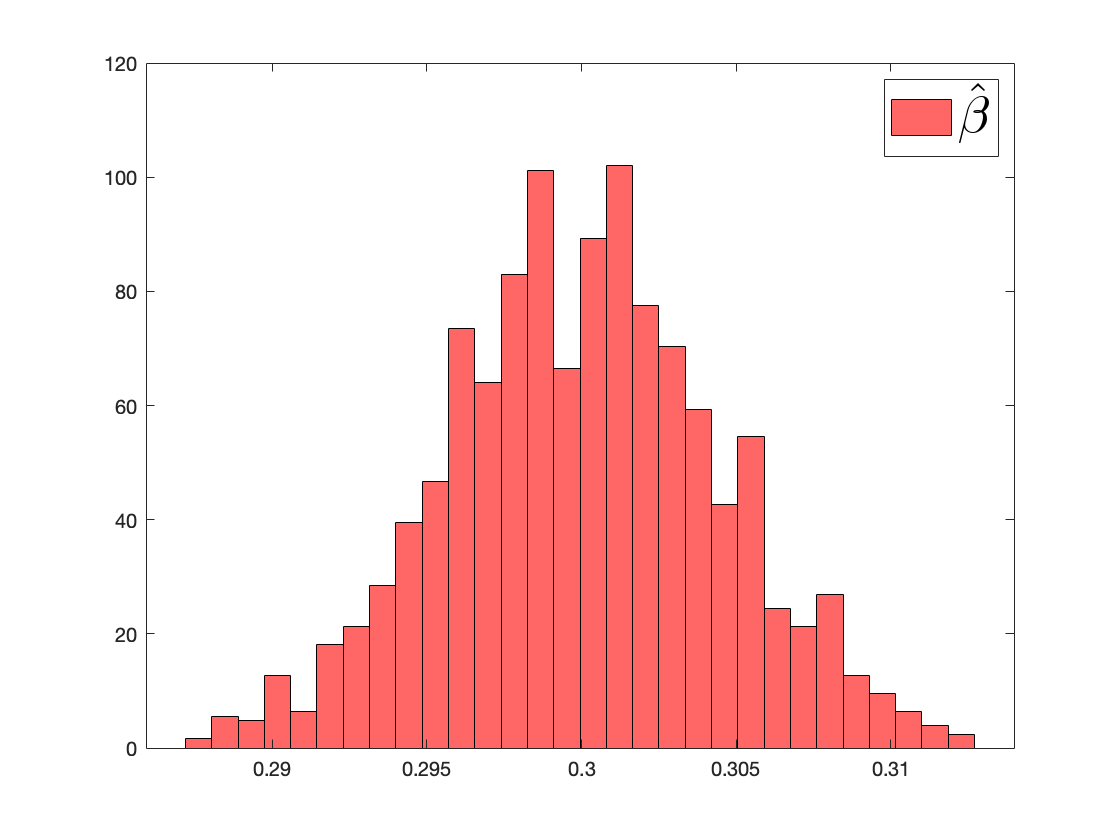}} 
    \subcaptionbox{$\hat{q}_2^{(L)} = 0.8001 \, (0.0285)$}
    {\includegraphics[width=0.49\linewidth]{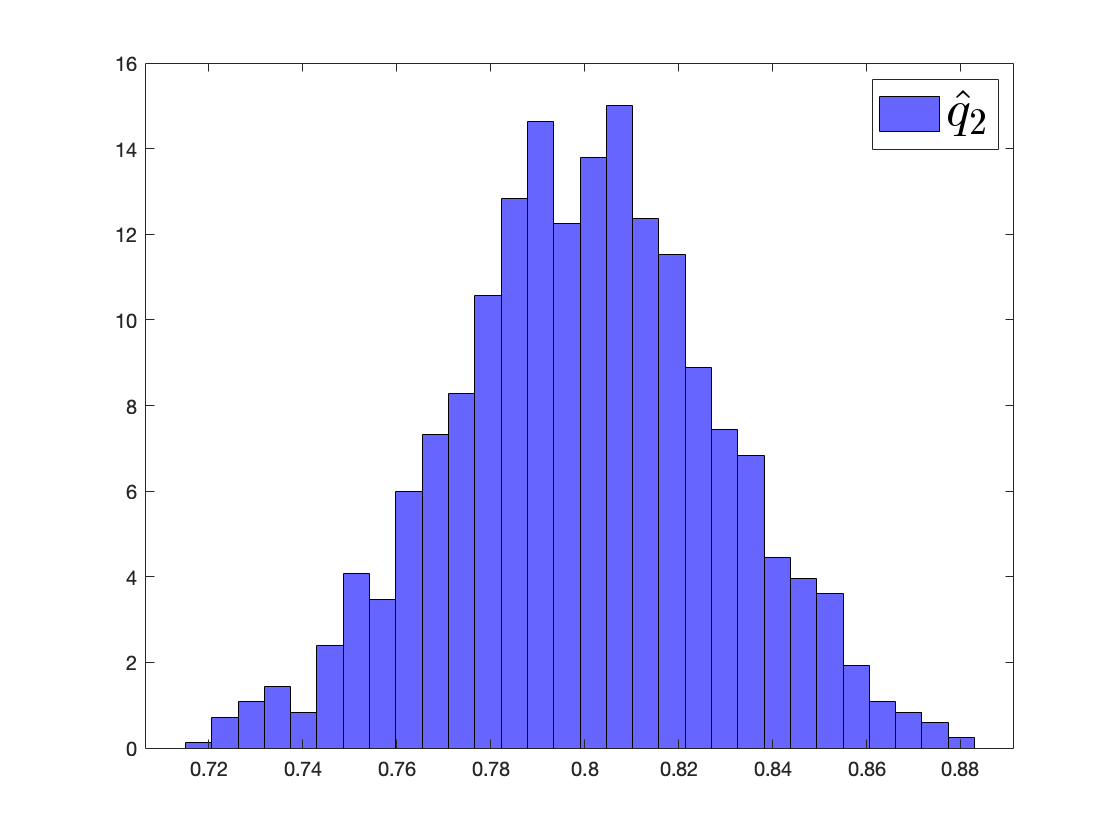}}
    \caption{Histograms of parameters in Case I.}
    \label{fig:sixTriStar}
\end{figure}

\begin{figure}
    \centering
    \subcaptionbox{$\hat{p}_0^{(L)} = 0.3000  \, (0.0019)$}
    {\includegraphics[width=0.49\linewidth]{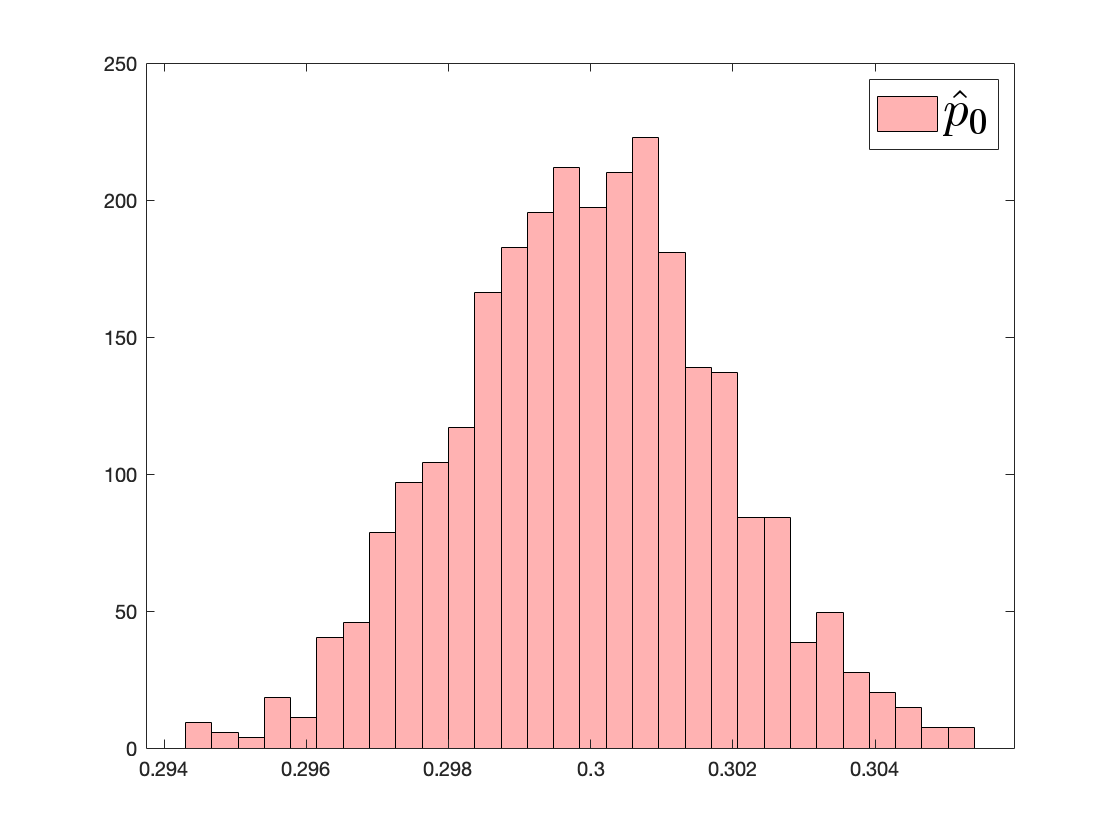}}
    \subcaptionbox{${\hat{q}_0^{(L)}} = 0.6001 \, (0.0028)$}
    {\includegraphics[width=0.49\linewidth]{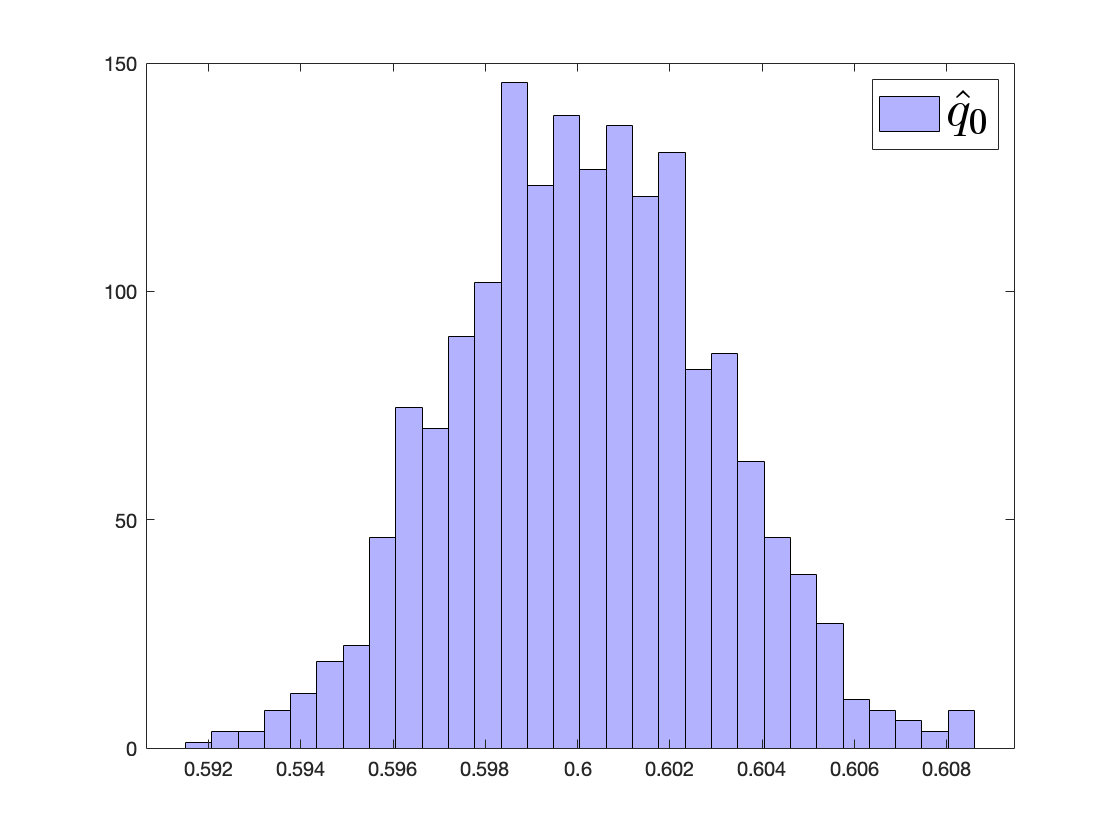}}\\
    \subcaptionbox{$\hat{\alpha}^{(L)} =  0.5004 \, (0.0114)$}
    {\includegraphics[width=0.49\linewidth]{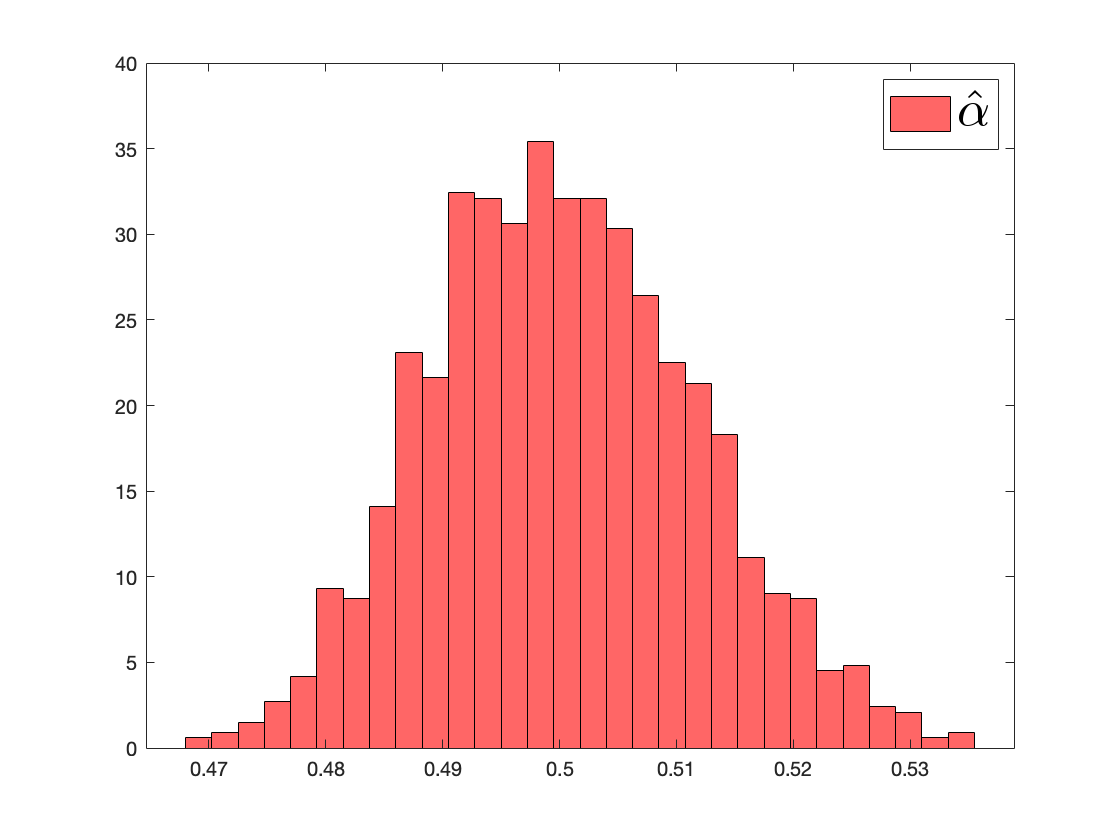}} 
    \hfill
    \subcaptionbox{$\hat{q}_1^{(L)} = 0.4000 \, (0.0059)$}
    {\includegraphics[width=0.49\linewidth]{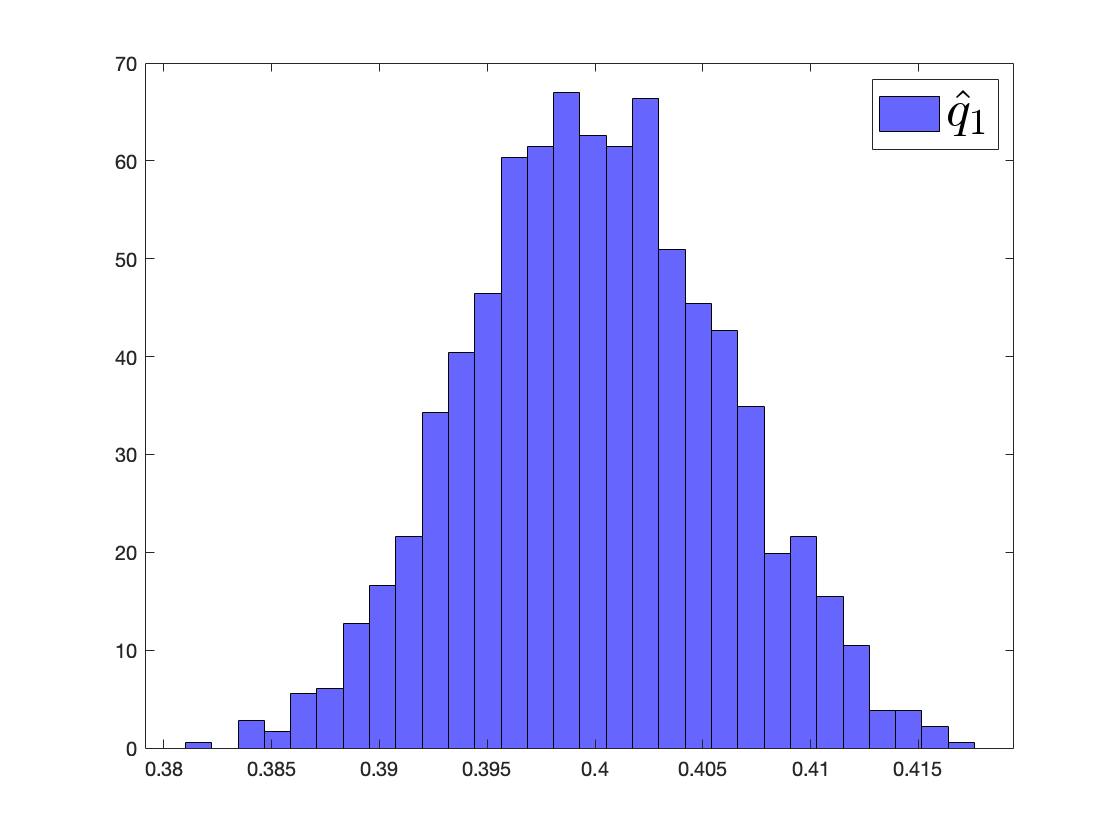}} \\
    \subcaptionbox{$\hat{\beta}^{(L)} = 0.3001 \, (0.0035)$}
    {\includegraphics[width=0.49\linewidth]{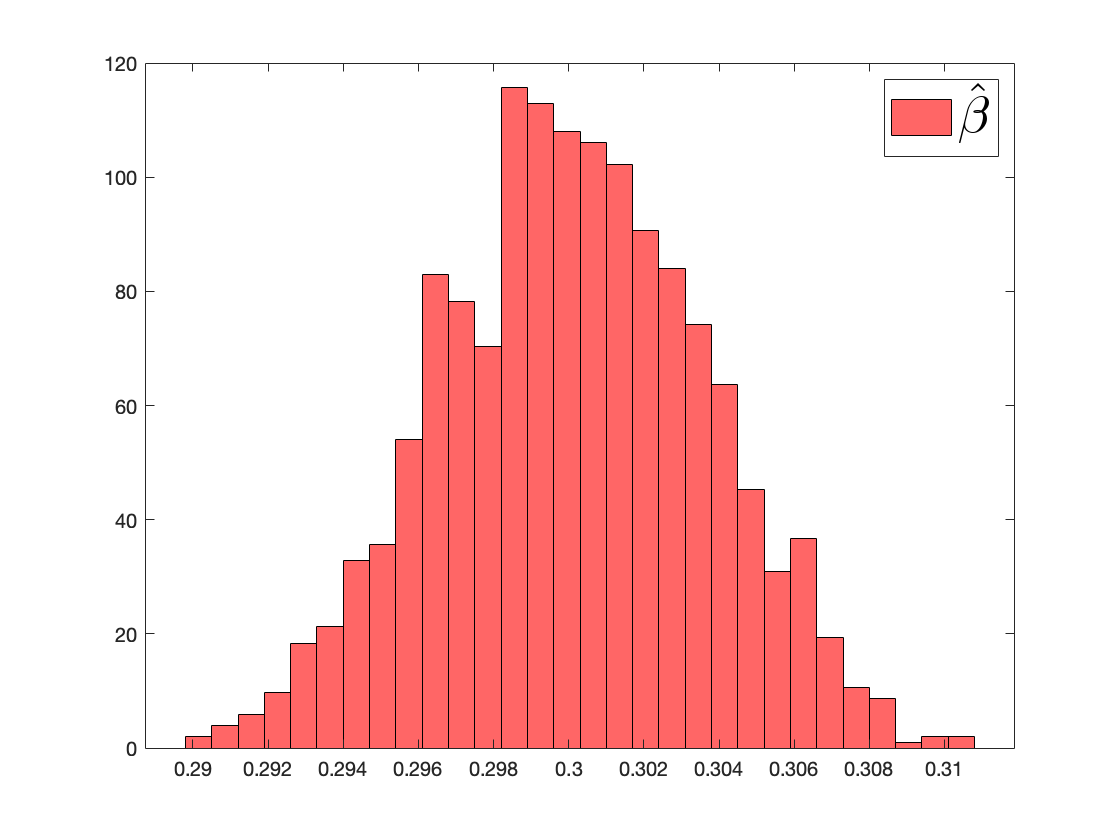}} 
    \subcaptionbox{$\hat{q}_2^{(L)} = 0.8006 \, (0.0226)$}
    {\includegraphics[width=0.49\linewidth]{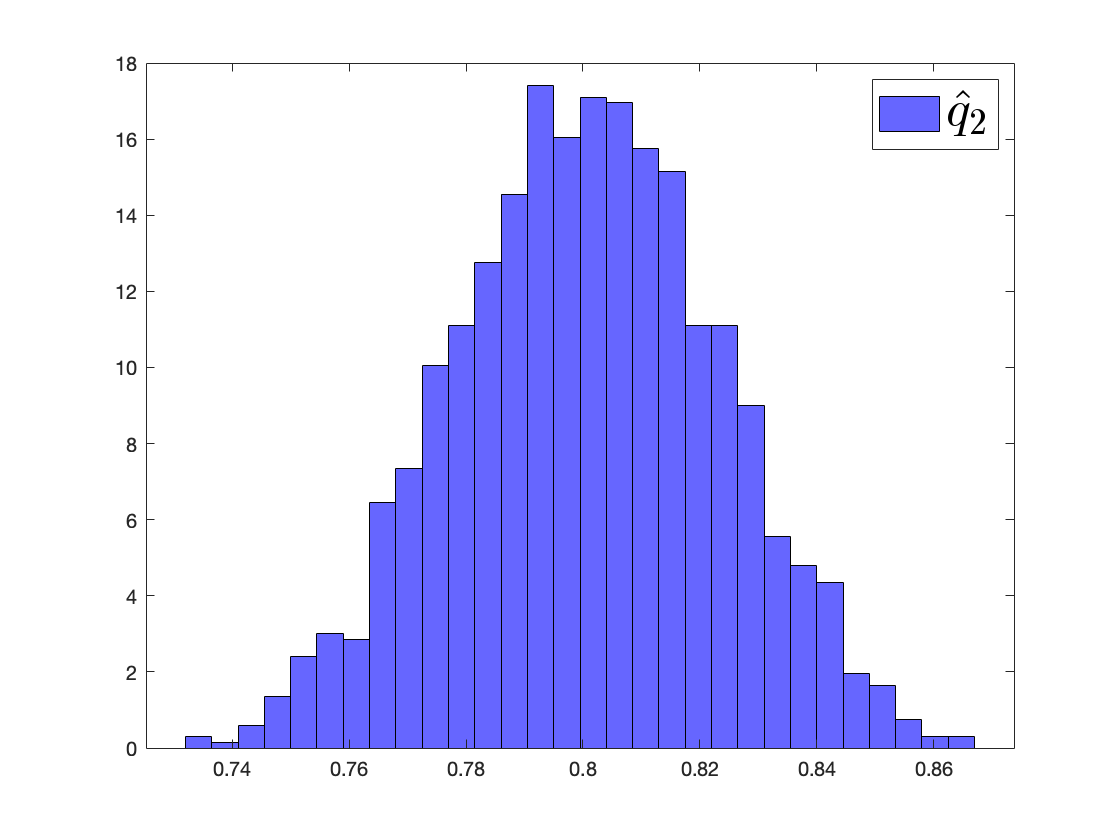}}
    \caption{Histograms of parameters in Case II.}
    \label{fig:sixStarStar}
\end{figure}

\begin{figure}
    \centering
    \subcaptionbox{$\hat{\alpha}^{(L)} = 0.5322 \, (0.1660)$}
    {\includegraphics[width=0.49\linewidth]{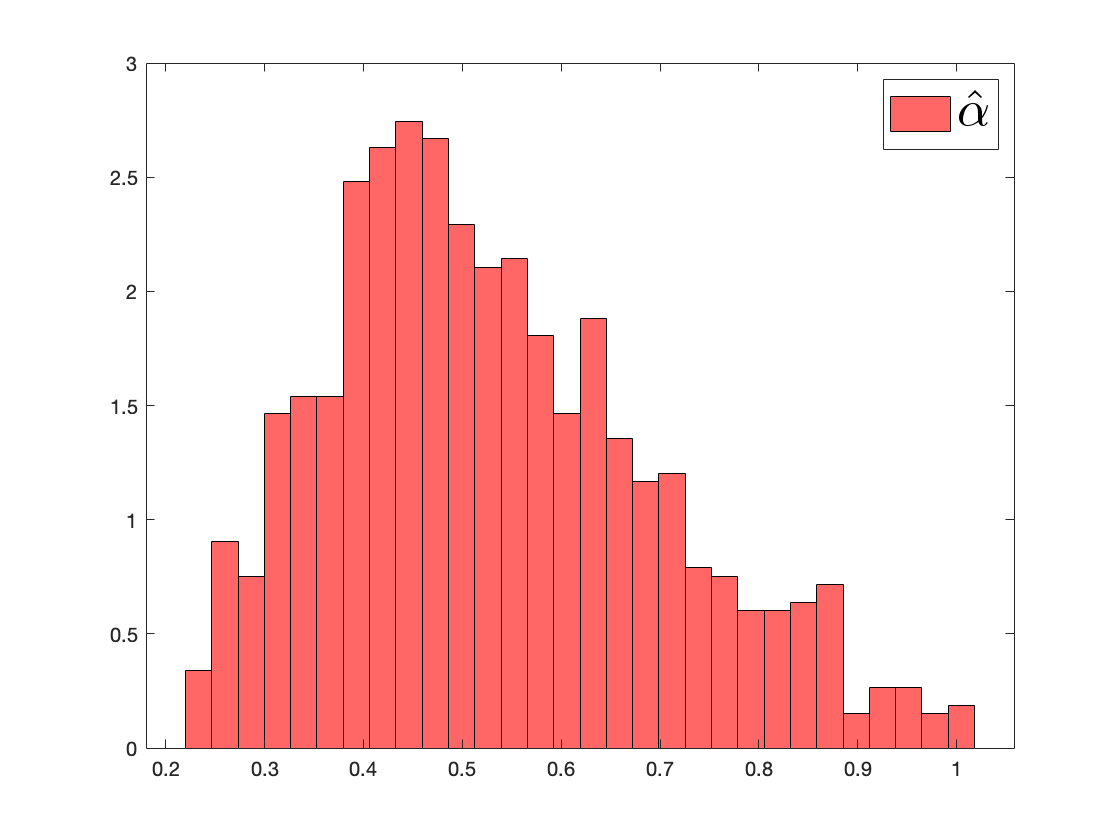}}
    \hfill
    \subcaptionbox{$\hat{\lambda}^{(L)} = 1.5232 \, (0.3238)$}
    {\includegraphics[width=0.49\linewidth]{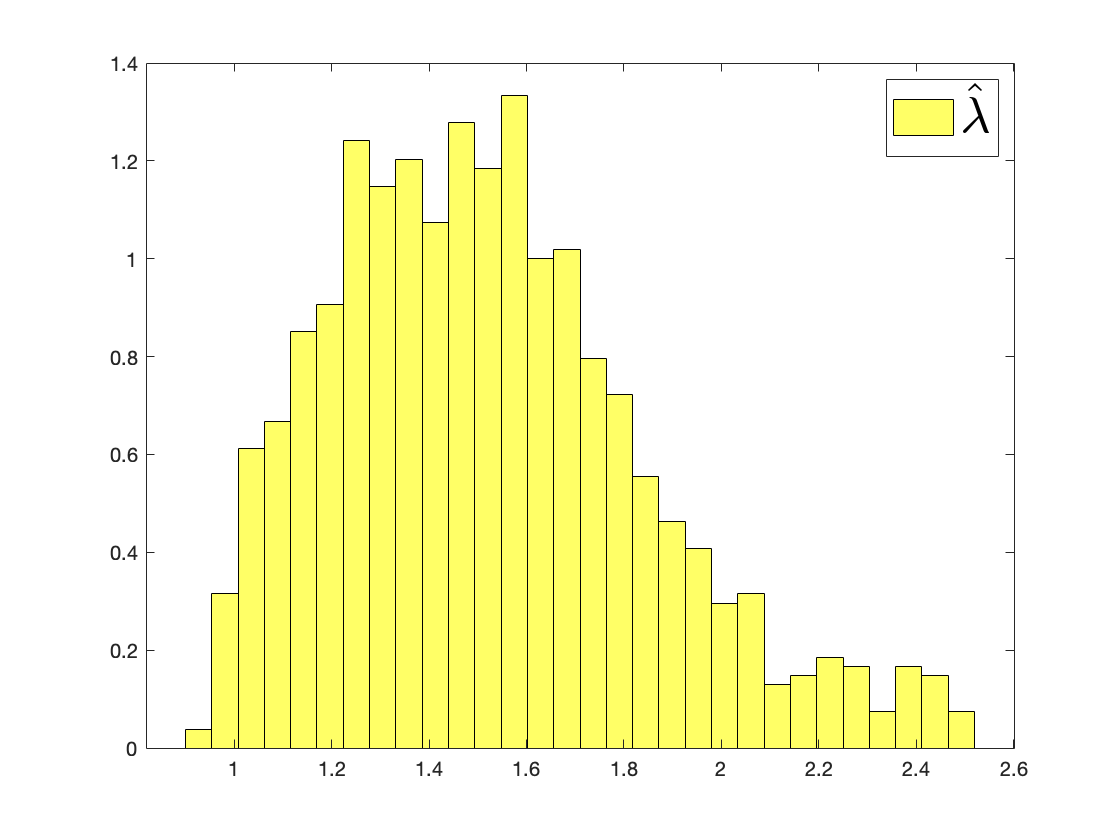}} 
    \caption{Histograms of parameters in Case III}
    \label{fig:seven}
\end{figure}

For every parameter instance considered, we perform $R$ runs, each run corresponding to $T$ time epochs. Let $\hat{p}_0^{\,r}$ denote the estimate produced in the $r$-th run, with $r \in \{1,2,\ldots,R\}$. We in addition define
\[
\Bar{p}_0^{(R)} := \frac{\sum_{r = 1}^R\hat{p}_0^{\,r}}{R} \qquad \Bar \sigma_0^{(R)} := \sqrt{\frac{\sum_{r = 1}^R (\hat{p}_0^{\,r} - \Bar{p}_0^{(R)})^2}{R}}
\]
as the sample mean and sample standard deviation of the estimates of $p_0$ that result from the $R$ runs, respectively. Similar notations have been used in relation to the other parameters. In our experiments, we systematically use $R = 1\,000$ and $T = 10^5$. The number of vertices is $N = 15$, so that the number of possible edges is $n = 105$. Importantly, all estimates are based on the same simulated trace, so as to make the comparison maximally fair. In Figures \ref{fig:sixTriStar}, \ref{fig:sixStarStar}, and \ref{fig:seven}, histograms of the $R$ estimates are presented; they cover each of the parameters in the three cases considered. {In this context, it is noted that the estimates $\hat p_0$, $\hat q_0$, $\hat q_1$, $\hat \beta$, $\hat q_2$ are {\it identical} for Cases I and III;
this property is explained in detail below. Therefore, in Figure \ref{fig:seven} (pertaining to Case III) we only present the estimates of $\lambda$ and $\alpha$.}
The sample mean and (between brackets) the sample standard deviation are given below each of the histograms. To evaluate the scalability of our method, we also consider the case $N = 30$, and present the corresponding values of $\Bar{p}_0^{(R)}$ and $\Bar \sigma_0^{(R)}$.

\medskip

Each run consists of a simulation phase followed by an estimation phase, both executed within the MATLAB environment. In each run, we simulate the process once, record the occurrences of $K_2$, $K_3$, $S_3$, and $S_4$, and compute both the single-snapshot moments and the cross-temporal moments. We then apply the two-step estimation procedure described in Section~\ref{ssec:est} to estimate the corresponding unknown parameters. The simulation phase takes approximately 44 seconds, while the estimation phase never exceeds 5 seconds.

General conclusions from Figures \ref{fig:sixTriStar}, \ref{fig:sixStarStar}, and \ref{fig:seven} ($N = 15$) are that (i)~all parameters are reproduced with a relatively high level of precision, and (ii)~the bell-shaped curves suggest that the estimators may be approximately Gaussian. For comparison, we also present the estimates for $N = 15$ and $N = 30$ in Table~\ref{tab:compare}.
\begin{table}[h!]
\centering
\renewcommand{\arraystretch}{1.5}
\begin{tabular}{|c|c||c|c|}
\hline
\multicolumn{2}{|c||}{\textbf{Statistic}} & \textbf{$N=15$} & \textbf{$N=30$} \\
\hline \hline
\multirow{6}{*}{Case I} 
& $\hat{p}_0^{(L)}$ & $0.2999  \, (0.0020)$ & $0.3000 \, (0.0017)$ \\
& $\hat{q}_0^{(L)}$ & $0.6002 \, (0.0030)$ & $0.6000 \, (0.0027)$\\
& $\hat{\alpha}^{(L)}$ & $0.5003 \,(0.0127)$ & $0.5014 \, (0.0214)$\\
& $\hat{q}_1^{(L)}$ & $0.4000 \, (0.0059)$ & $0.4001 \, (0.0111)$\\
& $\hat{\beta}^{(L)}$ & $0.3000 \,(0.0044)$ & $0.3005, \, (0.0070)$ \\
& $\hat{q}_2^{(L)}$ & $0.8001\, (0.0285)$ & $0.8028 \, (0.0444)$ \\
\hline
\multirow{6}{*}{Case II} & $\hat{p}_0^{(L)}$ & $0.3000  \, (0.0019)$ & $0.2999 \, (0.0017)$\\
& $\hat{q}_0^{(L)}$ & $0.6001 \, (0.0028)$ & $0.6000 \, (0.0027)$ \\
& $\hat{\alpha}^{(L)}$ & $0.5004 \, (0.0114)$ & $0.5013 \, (0.0203)$ \\
& $\hat{q}_1^{(L)}$ & $0.4000 \, (0.0059)$ & $0.4001 \, (0.0104)$\\
& $\hat{\beta}^{(L)}$ & $0.3001 \, (0.0035)$ & $0.3004 \, (0.0060)$\\
& $\hat{q}_2^{(L)}$ & $0.8006 \, (0.0226)$ & $0.8022, \, (0.0381)$ \\
\hline
\hline
\multirow{2}{*}{Case III} & $\hat{\alpha}^{(L)}$ & $0.5322 \, (0.1660)$ & $0.8304 \, (0.5848)$ \\
& $\hat{\lambda}^{(L)}$ & $1.5232 \, (0.3238)$ & $1.4089 \, (0.6460)$ \\
\hline
\end{tabular}
\caption{Comparison of estimates for $N = 15, \, 30$. }
\label{tab:compare}
\end{table}

We continue by presenting more detailed observations, starting with a comparison between Case~I and Case II. 
It is seen that the estimates obtained for Case II have smaller standard deviations than their counterparts for Case~I. In Case I, we use triangles ($K_3$), while in Case II, we use 4-stars (i.e., $S_4$). They both require three edges, but the number of triangles is, under our parameter setting, substantially lower than 4-stars; it may be this low number of triangles in Case I that makes the estimator perform less well. (For example, in a single run, the mean  number of $K_3$ is about 75, whereas the mean number of $S_4$ is about 900.)

Comparing Case I with Case III, we find that the stationary probabilities $\hat \pi_1$, $\hat \varrho_1$, and $\hat \varrho_2$ are the same in both cases, as they are derived using identical moment equations in Step 1. In Step 2 of Case I, there are three equations involving $\mathcal{A}_{T,1}$, $\mathcal{K}_{T,3,1}$, $\mathcal{S}_{T,3,1}$ with three unknowns: 
\[
{\mathbb{E} \, Z_1}, \qquad {\mathbb{E} \, X_1}, \qquad {\mathbb{E} \, Y_1} \,.
\]
In Case III, to match the number of equations and the number of unknown parameters, we use an additional moment equation involving $\mathcal{A}_{T,2}$. The unknowns in this case are:
\[
{\mathbb{E} \, Z_1}, \qquad {\mathbb{E} \, X_1}, \qquad {\mathbb{E} \, Y_1}, \qquad 
\bar f_1(2)\,.
\]
This means that the estimates of 
${\mathbb{E} \, Z_1}$, 
${\mathbb{E} \, X_1}$, 
${\mathbb{E} \, Y_1}$ are the same in Cases I and III, as they are obtained by solving the same set of equations. Consequently, ${\mathbb{E} \, Z_2}$, 
${\mathbb{E} \, X_2}$, 
${\mathbb{E} \, Y_2}$ are also identical, as they are derived from
${\mathbb{E} \, Z_1}$, 
${\mathbb{E} \, X_1}$, 
${\mathbb{E} \, Y_1}$, along with 
$\hat \pi_1$, $\hat \varrho_1$, $\hat \varrho_2$. It follows that the five parameters $\hat p_0$, $\hat q_0$, $\hat q_1$, $\hat \beta$, $\hat q_2$ are the same in both cases, as they are computed using ${\mathbb{E} \, Z_1}$, ${\mathbb{E} \, Z_2}$, ${\mathbb{E} \, Y_1}$, ${\mathbb{E} \, X_2}$, ${\mathbb{E} \, Y_2}$, respectively.
In Case I, $\hat \alpha$ can be directly calculated from ${\mathbb{E} \, X_1}$. However, in Case III, ${\mathbb{E} \, X_1}$ depends on both $\lambda$ and $\alpha$.  Since $\bar f_1(2) = e^{-\lambda}/{\mathbb{E} \, X_1}$, $\lambda$ and $\alpha$ must be estimated jointly from $\bar f_1(2)$ and ${\mathbb{E} \, X_1}$. This additional complexity explains the substantially larger standard  deviation observed in the estimate of $\alpha$.

\section{Conclusion}
\label{sec:conclusion}

In this paper, we have studied a modulated version of the dynamic Erd\H{o}s-R\'enyi random graph. Concretely, at any point in time the mode of an underlying background process determines which of two dynamic Erd\H{o}s-R\'enyi random graphs becomes visible. 
Our goal is to estimate the parameters of the two dynamic Erd\H{o}s-R\'enyi random graphs and the background process, solely based on subgraph counts pertaining to the observed process.
Our estimation approach is based on the method of moments. The main idea is that we find closed-form expressions for various moments (involving the observations of the subgraph counts over time), equate them to their estimated counterparts, and solve for the unknown parameters.

\medskip

Potential directions for future research include:
\begin{itemize} \item[$\circ$] 
Is it possible to identify subgraphs which consistently provide `more information' than others, i.e., consistently lead to estimators with lower variance?
Some subgraphs are relatively rare (for instance $\ell$-complete graphs for higher values of $\ell$), and therefore less suitable for estimation purposes. 
It is noted that it is in this respect also a consideration that some subgraphs are easier to count (from the graph's adjacency matrix, that is) than others; for instance counting stars just requires knowing each vertex's number of neighbors. 

\item[$\circ$] 
 Besides selecting the set of subgraphs, a second set of decisions concerns the choice of the `lag' $d$. For higher values of $d$ the subgraph counts at times $k$ and $k+d$ are virtually independent, so that the corresponding cross moment provides a relatively low amount of information.

  \item[$\circ$] The number of moment equations must be at least the number of unknown parameters, and the equations need to be `as independent as possible'. Is it a good strategy to work with more moment equations than parameters, and then apply least squares to obtain the estimates? 
 \item[$\circ$] Can we prove properties of the resulting estimators? The histograms that we presented suggest the estimators to be asymptotically normal. \end{itemize}

\end{document}